\documentclass[reqno,twoside,11pt]{amsart}
\usepackage{cite}
\usepackage{amsmath}
\usepackage{amsfonts}
\usepackage{amssymb}
\usepackage{epsfig}
\usepackage{verbatim}
\usepackage{color}

\IfFileExists{epsf.def}{\input epsf.def}{\usepackage{epsf}}

\IfFileExists{mathptmx.sty}{\usepackage{mathptmx}}{\usepackage{mathpazo}}
\usepackage{mathrsfs}

\DeclareFontFamily{OT1}{eusb}{} \DeclareFontShape{OT1}{eusb}{m}{n} {<5> <6> <7> <8> <9> <10> <11> <12> <14.4> eusb10}{}
\DeclareMathAlphabet{\eusb}{OT1}{eusb}{m}{n}

\DeclareFontFamily{OT1}{eusm}{} \DeclareFontShape{OT1}{eusm}{m}{n} {<5> <6> <7> <8> <9> <10> <11> <12> <14.4> eusm10}{}
\DeclareMathAlphabet{\eusm}{OT1}{eusm}{m}{n}

\DeclareFontFamily{OT1}{eufm}{} \DeclareFontShape{OT1}{eufm}{m}{n} {<5> <6> <7> <8> <9> <10> <11> <12> <14.4> eufm10}{}
\DeclareMathAlphabet{\mathfrak}{OT1}{eufm}{m}{n}

\DeclareFontFamily{OT1}{fraktura}{}
\DeclareFontShape{OT1}{fraktura}{m}{n} {<5> <6> <7> <8> <9> <10> <11> <12> <13> <14.4> [1.1] eufm10}{}
\DeclareMathAlphabet{\fraktura}{OT1}{fraktura}{m}{n}

\DeclareFontFamily{OT1}{cmfi}{} \DeclareFontShape{OT1}{cmfi}{m}{n} {<5> <6> <7> <8> <9> <10> <11> <12> <13> <14.4> [0.9] cmfi10}{}
\DeclareMathAlphabet{\cmfi}{OT1}{cmfi}{b}{n}

\DeclareFontFamily{OT1}{cmss}{} \DeclareFontShape{OT1}{cmss}{m}{n} {<5> <6> <7> <8> <9> <10> <11> <12> <13> <14.4> cmss10}{}
\DeclareMathAlphabet{\cmss}{OT1}{cmss}{m}{n}

\newtheoremstyle{thm}{1.5ex}{1.5ex}{\itshape\rmfamily}{} {\bfseries\rmfamily}{}{2ex}{}

\newtheoremstyle{def}{1.5ex}{1.5ex}{\rmfamily\sl}{} {\bfseries\rmfamily}{}{2ex}{}

\newtheoremstyle{rem}{1.3ex}{1.3ex}{\rmfamily}{} {\bfseries\rmfamily}{}{2ex}{}

\newtheoremstyle{ass}{1.5ex}{1.5ex}{\rmfamily\sl}{} {\bfseries\rmfamily}{}{2ex}{}

\newenvironment{proofsect}[1] {\vskip0.1cm\noindent{\rmfamily\itshape#1.}}{\qed\vspace{0.15cm}}

\theoremstyle{thm}
\newtheorem{theorem}{Theorem}[section]
\newtheorem{lemma}[theorem]{Lemma}
\newtheorem{proposition}[theorem]{Proposition}

\newtheorem*{Main Theorem}{Main Theorem.}
\newtheorem{corollary}[theorem]{Corollary}

\newtheoremstyle{named}{}{}{\itshape}{}{\bfseries}{}{.5em}{\thmnote{#3}}
\theoremstyle{named}

\theoremstyle{def}

\theoremstyle{rem}
\newtheorem{remark}[theorem]{{Remark}}

\numberwithin{equation}{section}

\renewcommand{\section}{\secdef\sct\sect}
\newcommand{\sct}[2][default]{\refstepcounter{section}
\addcontentsline{toc}{section}
{{\tocsection {}{\thesection}{\!\!\!\!#1\dotfill}}{}}
\vspace{0.7cm}
\centerline{ 
\scshape\arabic{section}.\ #1} \nopagebreak \vspace{0.2cm}}
\newcommand{\sect}[1]{
\vspace{0.4cm} \centerline{\large\scshape\rmfamily #1}
\vspace{0.2cm}}

\renewcommand{\subsection}{\secdef\subsct\sbsect}
\newcommand{\subsct}[2][default]{\refstepcounter{subsection}
\addcontentsline{toc}{subsection}
{{\tocsection{\!\!}{\hspace{1.2em}\thesubsection}{\!\!\!\!#1\dotfill}}{}}
\nopagebreak\vspace{0.45\baselineskip} {\flushleft\bf
\arabic{section}.\arabic{subsection}~\bf #1.~}
\\*[3mm]\noindent
\nopagebreak}
\newcommand{\sbsect}[1]{\vspace{0.1cm}\noindent
\textbf{#1.~}\vspace{0.1cm}}

\renewcommand{\subsubsection}{%
\secdef \subsubsect\sbsbsect}
\newcommand{\subsubsect}[2][default]{%
\refstepcounter{subsubsection} 
\addcontentsline{toc}{subsubsection}{{\tocsection{\!\!}
{\hspace{3.05em}\thesubsubsection}{\!\!\!\!#1\dotfill}}{}}
\nopagebreak
\vspace{0.15\baselineskip} \nopagebreak {\flushleft\rmfamily
\itshape\arabic{section}.\arabic{subsection}.\arabic{subsubsection}
\ \rmfamily #1\/.}\ }
\newcommand{\sbsbsect}[1]{\vspace{0.1cm}\noindent
\rmfamily \itshape
\arabic{section}.\arabic{subsection}.\arabic{subsubsection} \
\sffamily #1\/.\ }

\renewcommand{\caption}[1]{%
\vglue0.5cm
\refstepcounter{figure}
\begin{minipage}{0.9\textwidth}\small {\sc Figure~\thefigure. }#1\end{minipage}}

\newcommand{\dist}{\operatorname{dist}}

\newcommand{\textd}{\text{\rm d}\mkern0.5mu}

\newcommand{\texte}{\text{\rm  e}\mkern0.7mu}

\newcommand{\Var}{\text{\rm Var}}
\newcommand{\Cov}{\text{\rm \Cov}}
\newcommand{\1}{{1\mkern-4.5mu\textrm{l}}}
\renewcommand{\1}{\text{\sf 1}}

\newcommand{\FF}{\mathcal F}

\newcommand{\NN}{\mathcal N}

\newcommand{\D}{\mathbb D}

\newcommand{\N}{\mathbb N}

\newcommand{\Q}{\mathbb Q}
\newcommand{\R}{\mathbb R}

\newcommand{\Z}{\mathbb Z}

\newcommand{\twoeqref}[2]{(\ref{#1}--\ref{#2})}

\newcommand{\cc}{{\text{\rm c}}}

\newcommand{\Lawarrow}{{\,\overset{\text{\rm law}}\longrightarrow\,}}
\newcommand{\Lawlongarrow}{{\,\,\overset{\text{\rm law}}\longrightarrow\,\,}}

\def\myffrac#1#2 in #3{\raise 2.6pt\hbox{$#3 #1$}\mkern-1.5mu\raise 0.8pt\hbox{$#3/$}\mkern-1.1mu\lower 1.5pt\hbox{$#3 #2$}}
\newcommand{\ffrac}[2]{\mathchoice%
	{\myffrac{#1}{#2} in \scriptstyle}
	{\myffrac{#1}{#2} in \scriptstyle}
	{\myffrac{#1}{#2} in \scriptscriptstyle}
	{\myffrac{#1}{#2} in \scriptscriptstyle}
}

\newcommand{\wt}{\widetilde}

\newcommand{\laweq}{\,\overset{\text{\rm law}}=\,}
\newcommand{\argmax}{{\text{\rm argmax}}\,}

\newcommand{\hf}{h^{\text{\rm f}}}
\newcommand{\hc}{h^{\text{\rm c}}}

\begin{document}

\title[Extreme points of 2D Gaussian free field\hfill]
{\large Extreme local extrema of two-dimensional\\discrete Gaussian free field}
\author[\hfill M.~Biskup and O.~Louidor]
{Marek~Biskup$^1$ and\, Oren~Louidor$^{1,2}$}
\thanks{\hglue-4.5mm\fontsize{9.6}{9.6}\selectfont\copyright\,\textrm{2015} \textrm{M.~Biskup, O.~Louidor.
Reproduction, by any means, of the entire
article for non-commercial purposes is permitted without charge.\vspace{2mm}}}
\maketitle

\vspace{-5mm}
\centerline{\textit{$^1$Department of Mathematics, UCLA, Los Angeles, California, USA}}
\centerline{\textit{$^2$Faculty of Industrial Engineering and Management, Technion, Haifa, Israel}}

\vspace{3mm}
\begin{quote}
\footnotesize \textbf{Abstract:}
We consider the discrete Gaussian Free Field in a square box in~$\Z^2$ of side length~$N$ with zero boundary conditions and study the joint law of its properly-centered extreme values~($h$) and their scaled spatial positions ($x$) in the limit as $N\to\infty$. Restricting attention to extreme local maxima, i.e., the extreme points that are maximal in an $r_N$-neighborhood thereof, we prove that the associated process tends, whenever $r_N\to\infty$ and $r_N/N\to0$, to a Poisson point process with intensity measure $Z(\textd x)\texte^{-\alpha h}\textd h$, where $\alpha:=
2/\sqrt{g}$ with $g:=2/\pi$ and where~$Z(\textd x)$ is a random Borel measure on $[0,1]^2$. In particular, this yields an integral representation of the law of the absolute maximum, similar to that found in the context of Branching Brownian Motion. We give evidence that the random measure~$Z$ is a version of the derivative martingale associated with the continuum Gaussian Free Field. 
\end{quote}
\vspace{-3mm}

\section{Introduction}
\vspace{-4mm}
\subsection{Main results}
\noindent
Consider a box~$V_N:=(0,N)^2\cap\Z^2$ in the square lattice and let $G_N(x,y)$ denote the Green function of the simple symmetric random walk started from~$x$ and killed upon exiting~$V_N$. 
The two-dimen\-sio\-nal Discrete Gaussian Free Field (DGFF) in~$V_N$ is a collection of Gaussian random variables $\{h_x\colon x\in V_N\}$ with mean zero and covariance $\text{Cov}(h_x,h_y):=G_N(x,y)$. Another way to define the DGFF is by prescribing its full distribution; this is achieved by normalizing the measure
\begin{equation}
\label{E:1.1}
\exp\Bigl\{-\frac18\sum_{\langle x,y\rangle}(h_x-h_y)^2\Bigr\}\prod_{x\in V_N}\textd h_x\prod_{x\in\partial V_N}\delta_0(\textd h_x).
\end{equation}
Here the sum goes over unordered nearest-neighbor pairs with at least one vertex in~$V_N$ and the product of Dirac delta's imposes a Dirichlet boundary condition on the outer boundary $\partial V_N$ of~$V_N$. By \eqref{E:1.1} the DGFF has the \emph{Gibbs-Markov property}: Conditional on $\{h_z\colon z\ne x\}$, the field~$h_x$ reduced by the average of~$h_z$ over the nearest neighbors~$z$ of~$x$ has the law of a standard normal.

The aim of this paper is to study the statistics of extreme values of the DGFF in the limit \hbox{$N\to\infty$}. We will focus attention on large local maxima, i.e., those extreme points whose value dominates the configuration in an $r$-neighborhood thereof. Thus, for $r\ge1$, let $\Lambda_r(x):=\{z\in\Z^2\colon |z-x|_1\le r\}$ and define a measure on $[0,1]^2\times\R$ by
\begin{equation}
\eta_{N,r}(A\times B):=\sum_{x\in V_N}\1_{\{x/N\in A\}}\1_{\{h_x-m_N\in B\}}\1_{\{h_x=\max_{z\in \Lambda_r(x)}h_z\}},
\end{equation}
for  Borel sets $A\subset[0,1]^2$ and $B\subset\R$ and a suitable centering sequence~$m_N$. A sample of~$\eta_{N,r}$ is a 
Radon measure supported on
a collection of points of the form $(x,h)$, where $x$ is the scaled location and~$h$ is the reduced height of a large ``peak'' in the underlying field configuration.

To study distributional limits, we endow the space of point measures on $[0,1]^2\times\R$ with the topology of vague convergence. For the centering sequence~$m_N$ we will take
\begin{equation}
\label{E:mN}
m_N:=2\sqrt g\,\log N-\frac34\sqrt g\,\log\log N,
\end{equation}
where $g:=2/\pi$ links $m_N$ to the asymptotic growth of the Green function which for~$x$ deep inside~$V_N$ scales as $G_N(x,x)=g\log N+O(1)$.
Anticipating Poisson limit laws, let us write ${\rm PPP}(\lambda)$ for the Poisson point process on a Polish space~$\Omega$ with sigma-finite intensity measure~$\lambda$. We will use this notation even when $\lambda$ is itself random (i.e., when ${\rm PPP}(\lambda)$ is a Cox process); the law of the points is then averaged over the law of~$\lambda$. Our principal result is then:

\begin{theorem}
\label{thm-main}
There is a random measure $Z(\textd x)$ on~$[0,1]^2$ with $Z([0,1]^2)<\infty$ a.s.\ and $Z(A)>0$ a.s.\ for any open set $A\subset[0,1]^2$ such that for any $r_N$ with $r_N\to\infty$ and $r_N/N\to0$, 
\begin{equation}
\label{E:1.9}
\eta_{N,r_N}\,\underset{N\to\infty}\Lawlongarrow\,\text{\rm PPP}\bigl(Z(\textd x)\otimes\texte^{-\alpha h}\textd h\bigr),
\end{equation}
where $\alpha:=2/\sqrt g$ --- which in present normalization reads $\alpha=\sqrt{2\pi}$. 
\end{theorem}

As an immediate consequence, we get information about the joint law of the (a.s.~unique) position and height of the absolute maximum:

\begin{corollary}
\label{cor1.5}
Let $\nu_N$ denote the law of $(N^{-1}\text{\rm argmax\,}h,\max_{x\in V_N}h_x-m_N)$ on~$[0,1]^2\times\R$.
For the random measure $Z(\textd x)$ from Theorem~\ref{thm-main}, define
\begin{equation}
\label{E:1.13b}
\widehat Z(A):=\frac{Z(A)}{Z([0,1]^2)}.
\end{equation}
Then $\nu_N\Lawarrow\nu$, where $\nu$ is for any Borel~$A\subset[0,1]^2$ given by
\begin{equation}
\label{E:1.14q}
\nu\bigl(A\times(-\infty,t]\bigr):= E\bigl(\,\widehat Z(A)\texte^{-\alpha^{-1}Z\texte^{-\alpha t}}\bigr),\qquad t\in\R \,,
\end{equation}
and $Z := Z([0,1]^2)$.
\end{corollary}

From \eqref{E:1.14q} we get an integral representation for the limit law of the centered maximum 
\begin{equation}
\label{E:1.9b}
P\bigl(\,\max_{x\in V_N}h_x-m_N\le t\bigr)\,\underset{N\to\infty}\longrightarrow\, E\bigl(\texte^{-\alpha^{-1}Z\texte^{-\alpha t}}\bigr),\qquad t\in\R,
\end{equation}
where $Z:=Z([0,1]^2)$. Letting $t\to\infty$ in \eqref{E:1.14q} in turn tells us that the limit law of $N^{-1}\argmax h$ is given by the probability measure $E(\widehat Z(\cdot))$.

Writing $Z(\textd x)\texte^{-\alpha h}\textd h$ as the product of~$\widehat Z(\textd x)$ and $Z\texte^{-\alpha h}\textd h$, where $Z:=Z([0,1]^2)$, the expression in \eqref{E:1.14q} suggests a \emph{sampling method} for the associated Poisson point process in \eqref{E:1.9}: First choose the heights from the Gumbel process with intensity $Z\texte^{-\alpha h}\textd h$ and then assign the spatial coordinates independently from $\widehat Z$. The spatial location of all large local maxima are thus governed by~$\widehat Z$. 
Of course,~$\widehat Z(A)$ is correlated with~$Z$ so explicit information is hard to get.

We also note that the measure~$Z(\textd x)$ is completely determined by (i.e., it is a measurable function of) a.e.\ sample of the limiting process~$\eta$. Indeed, if~$A\subset[0,1]^2$ is a Borel set and $(x_i,h_i)$ enumerates the points in $\eta$ with $x_i\in A$ with~$h_i$ listed in decreasing order, then
\begin{equation}
\alpha\,n\,\texte^{\alpha h_n}\,\underset{n\to\infty}\longrightarrow\, Z(A),\qquad\text{\rm a.s.}
\end{equation}
In particular, when only the field values (and not the positions) are of concern, the limiting law is that of a Gumbel process with intensity~$\texte^{-\alpha h}\textd h$, where all points are shifted by the random quantity $\alpha^{-1}\log Z([0,1]^2)$.

Finally, we remark that Theorem~\ref{thm-main} is only a first step towards the understanding of the full extreme point process associated with the DGFF. Indeed, each extreme local maximum will be surrounded by a ``cluster'' of points where the value of the field is only $O(1)$-term less than the value at the maximum. Naturally, such a situation should be described by a cluster process and handled by methods of two-scale convergence. Details of these will appear in a sequel paper. \textit{Update}: The convergence of the full cluster process has now been established and is the main subject of the forthcoming article~\cite{BL3}.

\subsection{The random measure}
\label{sec1.2}\noindent
One of the most interesting consequences of Theorem~\ref{thm-main} is the existence of the random measure~$Z(\textd x)$. Here are some facts about the statistical properties of this object: 

\begin{theorem}
\label{thm1.3}
The random measure $Z(\textd x)$ from Theorem~\ref{thm-main} is concentrated on~$(0,1)^2$ a.s., it is a.s.\ non-atomic and its total mass,~$Z:=Z([0,1]^2)$, satisfies the moment bounds
\begin{equation}
E(Z^p)\begin{cases}
<\infty,\qquad&\text{if }p\in(-\infty,1),
\\
=\infty,\qquad&\text{if }p\in[1,\infty).
\end{cases}
\end{equation}
For the borderline case $p:=1$ we in fact get that $Z(1\vee\log Z)^q$ is integrable for~$q<-1$ and not integrable for~$q\ge-1$.
\end{theorem}

\noindent
The positive moments can be controlled thanks to our knowledge of the behavior of the Laplace transform of~$Z:=Z([0,1]^2)$ for small values of its arguments:
\begin{equation}
\label{E:1.10}
E(\texte^{-\lambda Z})=1-[C_\star+o(1)]\lambda\log(\ffrac1 \lambda),\qquad \lambda\downarrow0,
\end{equation}
for some constant~$C_\star\in(0,\infty)$. Similar, albeit not so sharp, bounds exist also for the $\lambda\to\infty$ tail; these are still sufficient to control the negative moments of~$Z$ as well. (One in fact shows that $Z^{-1}$ has stretched-exponential moments.)

The measure~$Z$ can be constructed by a limiting procedure that uses objects that can more or less be explicitly identified. Let us begin by introducing proper notation. Given~$K\in\N$, consider the collection of $K^2$ disjoint boxes $B^{K,i}:=w_i^K+(0,\ffrac1K)^2$, where $w_1^K,\dots,w_{K^2}^K$ enumerates the points in $\{x/K\colon x\in \{0,\dots, K-1\}^2\}$. The closures of these boxes cover the closure of $B:=(0,1)^2$. Let $\Pi(x,A)$ denote the probability that the standard Brownian motion started at~$x\in B$ exits~$B$ through the set $A\subset\partial B$ and let $\Pi_i(x,A)$ denote the corresponding object for $B^{K,i}$, with $x\in B^{K,i}$ and $A\subset\partial B^{K,i}$. (Obviously, $\Pi$ and $\Pi_i$ are boundary harmonic measures in the respective sets.) For $x\in\bigcup_{i=1}^{K^2}B^{K,i}$, let $i(x)$ denote the unique index such that $x\in B^{K,i(x)}$. 

Writing $G(x,y)$ for the continuum Green function for the Laplacian on $(0,1)^2$ with zero boundary conditions, let $C_K\colon\bigcup_{i=1}^{K^2}B^{K,i}\times \bigcup_{i=1}^{K^2}B^{K,i}\to\R$ denote its harmonic extension, in both coordinates, to $[0,1]^2\setminus\bigcup_{i=1}^{K^2}B^{K,i}$. Explicitly,
\begin{equation}
C_K(x,y):=g\begin{cases}
\int_{\partial B}\Pi(x,\textd z)\log\frac{|z-y|}{|x-y|},\quad&\text{if }i(x)\ne i(y),
\\
\int_{\partial B}\Pi(x,\textd z)\log|z-y|-\int_{\partial B^{K,i}}\Pi_i(x,\textd z)\log|z-y|,&\text{if }i(x)=i(y)=i.
\end{cases}
\end{equation}
This function is symmetric and positive semi-definite and hence there exists a Gaussian process $\{\Phi_K(x)\colon x\in\bigcup_{i=1}^{K^2}B^{K,i}\}$ with mean zero and $\text{Cov}(\Phi_K(x),\Phi_K(y))=C_K(x,y)$. A.e.\ sample of the field $\Phi_K$ is harmonic, and thus infinitely differentiable throughout each $B^{K,i}$. Endowing the space of finite Borel measures on $[0,1]^2$ with the topology of vague convergence, we then get:

\begin{theorem}
\label{thm1.7}
Consider the measure~$Z(\textd x)$ from Theorem~\ref{thm-main} and let $\alpha:=2/\sqrt{g}=\sqrt{2\pi}$. There is a probability density~$\psi$ on~$[0,1]^2$, a constant $c_\star\in(0,\infty)$ and a sequence $\epsilon_K\downarrow0$ such that the random measure $Z_K(\textd x)$ on $[0,1]^2$, defined for $F(s):=s\texte^{-\alpha s}\1_{[0,\infty)}(s)$ by
\begin{equation}
\label{E:1.12t}
Z_K(A):=c_\star\,\sum_{i=1}^{K^2}\int_{(\epsilon_K,1-\epsilon_K)^2}\!\!\!\!\textd z\,\, \psi(z)\, F\bigl(2\sqrt{g}\log K-\Phi_K(w_i^K+\tfrac zK)\bigr)\1_{A}\bigl(w_i^K+\tfrac zK\bigr),
\end{equation}
obeys $Z_K \Lawarrow Z$ as $K \to \infty$.
\end{theorem}

Thanks to the form of the function $F$, and the fact that $\text{Var}(\alpha\Phi_K(z))=\alpha 2\sqrt{g}\log K+O(1)$, the measure in \eqref{E:1.12t} takes the form of a \emph{derivative martingale} (see Sections~\ref{sec2.4} and~\ref{sec2.5w}). 

\subsection{Scaling and conformal invariance}
The measure~$Z$ exhibits a host of interesting properties. We will now list some of these along with a heuristic explanation; rigorous proofs will appear in sequel papers.

The first observation is that the above conclusions apply to sets other than just square boxes. For a bounded open set~$D\subset\R^2$, consider the scaled-up lattice version $D_N:=\{x\in\Z^2\colon x/N\in D\}$. The DGFF in~$D_N$ is defined similarly as for the squares using the Green function of the simple random walk killed upon exit from~$D_N$. Proceeding similarly as for squares, the law of the spatial coordinate of the corresponding limit point process --- with heights still reduced by~$m_N$ in~\eqref{E:mN} --- should  be governed by a random measure~$Z^D(\textd x)$ on~$D$. 
The sheer existence (and uniqueness) of this object then implies
\begin{equation}
\label{E:1.18}
Z^{\lambda D}(\lambda\textd x)\laweq \lambda^4\, Z^D(\textd x),\qquad \lambda>0 \,.
\end{equation}
This follows from the asymptotics $\alpha(m_{\lambda N}-m_N)=4\,\log \lambda +o(1)$  valid for any $\lambda \in\Q$ and~$N$'s such that $\lambda N$ is a positive integer.  

In light of the Gibbs-Markov property of the DGFF, $Z^D$ also behaves quite predictably under partitions of~$D$. Indeed, suppose that $D', D''$ are two disjoint open subsets of $D$ such that $\overline{D'} \cup \overline{D''} = \overline{D}$. Setting $\gamma := \partial D' \cup \partial D''$, we then have
\begin{equation}
\label{E:1.19}
Z^D(\textd x)\laweq \texte^{\alpha\Phi(x)}Z^{D'}(\textd x)+\texte^{\alpha\Phi(x)}Z^{D''}(\textd x).
\end{equation}
Here $\Phi(x)$ is a random Gaussian field on $D$ which is the harmonic extension of the values on~$\gamma
$ of the Continuum Gaussian Free Field (CGFF) in~$D$ (alternatively, the centered Gaussian field on~$D$ whose covariance function is the harmonic extension in both coordinates of the Green function on~$D$ given its values on $\gamma$) while $Z^{D'}$ and~$Z^{D''}$ --- which we regard as measures on~$D$ --- are independent of each other and of~$\Phi(x)$. The argument uses~$Z^D(\gamma)=0$ a.s. 

The relations \twoeqref{E:1.18}{E:1.19}, along with the fact that the CGFF is conformally invariant imply a transformation rule for~$Z^D$ under conformal bijections of the underlying domain~$D$. Indeed, let $f\colon D\to \wt{D}$ be such a bijection. Then
\begin{equation}
\label{E:1.20}
(Z^{\wt{D}}\circ f)\,(\textd  x)\laweq |f'(x)|^4\,Z^D(\textd x).
\end{equation}
The measure $Z^D$ on any bounded and simply connected~$D$ can thus be obtained from that on the unit disc,~$Z^{\D}$, by a conformal map. Moreover, the law of the measure $(1-|x|^2)^{-4}Z^{\D}(\textd x)$ on~$\D$ is invariant under the M\"obius transforms preserving~$\D$.

We summarize these observations by noting that \twoeqref{E:1.18}{E:1.19} make~$Z^D$ an example of a \emph{Gaussian multiplicative chaos}; see Section~\ref{sec2.4} for some pointers to the literature. Property \eqref{E:1.20} is in turn a direct manifestation of the \emph{conformal invariance} of the continuum Gaussian Free Field at the level of extreme points. We emphasize that the role of conformal invariance for extreme points has been rather unclear because the continuum field does not exist as a function (so its maximum cannot be defined).

\noindent
\textit{Update}: The above properties of the~$Z$ measure have now been established rigorously in a large class of underlying domains (i.e., not just squares); cf Biskup and Louidor~\cite{BL2}.
\smallskip

\section{Earlier and related work}
\label{sec:RelatedWork}
\noindent
To keep the description of our main results succinct, we have insofar refrained from making any connections to earlier (and other related) work. Here we will mend this by giving a proper account of the literature that bears upon the subject at hand. We emphasize that not all what is quoted below is actually used in the paper. The principal novel ideas of this work will be discussed in Section~\ref{sec-3}.

\subsection{Gaussian Free Field in statistical mechanics}
The DGFF (in any spatial dimension) has been a source of much attention in equilibrium statistical mechanics. It arises in models of interfaces and/or crystal deformations (the so called harmonic crystal), fluctuation theories at critical points, field theory, etc. Thanks to the Dynkin~\cite{Dynkin} isomorphism theory (discovered already by Symanzik~\cite{Symanzik}), the DGFF is also closely connected with the local time of the simple random walk. Consequently, it bears upon subjects such as cover times of random walks and Markov chains (e.g., Dembo, Peres, Rosen and Zeitouni~\cite{DPRZ}, Ding, Lee and Peres~\cite{Ding-Lee-Peres}, Ding~\cite{Ding-cover-time}), random interlacements and loop decompositions (e.g. Le Jan~\cite{LeJan}, Sznitman~\cite{Sznitman}, Rodriguez~\cite{Rodriguez}), etc.

The two-dimensional DGFF is particularly interesting because its continuum version (CGFF) is invariant under conformal transformations of the underlying domain. This offers a framework for analyzing scaling limits of certain critical models. For instance, the level sets of the DGFF on the triangular lattice can be linked with the Schramm-Loewner process SLE$_4$ (Schramm and Sheffield~\cite{Schramm-Sheffield}), the height function associated with domino tilings scales to the CGFF (Kenyon~\cite{Kenyon}) etc.  A complication associated with the CGFF is that, by its scale-invariant nature, the ``field'' exists only as a random distribution on an appropriate function space. 

\subsection{The maximum}
Let us now move to the subject of interest in the present paper, which is the behavior of the extreme values of the DGFF in the limit as~$N\to\infty$. A particular aspect of this, the maximum
\begin{equation}
M_N:=\max_{x\in V_N}h_x,
\end{equation}
has been studied very intensely. Indeed, Bolthausen, Deuschel and Giacomin~\cite{BDeuG} showed that $M_N/\log N\to2\sqrt g$. More recently, building on Bolthausen, Deuschel and Zeitouni~\cite{BDZ}, Bramson and Zeitouni~\cite{BZ} proved that the family $\{M_N-m_N\}_{N\ge1}$ is tight when~$m_N$ is as in \eqref{E:mN}. Finally, in a very recent development, Bramson, Ding and Zeitouni~\cite{BDingZ} proved the following result:

\begin{theorem}[Theorem~1.1 of \cite{BDingZ}]
\label{thm2.1}
$M_N-m_N$ converges in distribution to a non-degenerate random variable.
\end{theorem}

\noindent
\textit{Update}: As Ofer Zeitouni informed us, although not explicitly derived in~\cite{BDingZ}, the representation of the limiting distribution as a Laplace transform of a random variable $Z$, i.e., formula~\eqref{E:1.9b}, can be extracted from their work as well.
\smallskip

In two other recent papers, Ding~\cite{Ding} and Ding and Zeitouni~\cite{DZ} have studied the tails of the maximum. Specifically, for the upper tail they derived
\begin{equation}
\label{E:2.2q}
C^{-1}
t\texte^{-\sqrt{2\pi}\,t}\le P\bigl(M_N\ge m_N+t\bigr)\le C t\texte^{-\sqrt{2\pi}\,t},\qquad 1\le t < \sqrt{\log N},
\end{equation}
for some $C\ge1$. For the lower tail, they got the estimates
\begin{equation}
\label{E:2.3q}
c\texte^{-\texte^{Ct}}\le P\bigl(M_N\le m_N-t\bigr)\le C\texte^{-\texte^{ct}},\qquad 0\le t\le(\log N)^{2/3},
\end{equation}
for some $0<c<C<\infty$. The recent work of Bramson, Ding and Zeitouni~\cite{BDingZ} controls the asymptotic form of the upper tail of~$M_N$ including the multiplicative constant. Indeed, recasting Proposition~2.2 of~\cite{BDingZ} into a simpler form, we get:

\begin{theorem}[Proposition~2.2 of~\cite{BDingZ}]
\label{thm-asymp}
There exists a constant $C_\star\in(0,\infty)$ such that
\begin{equation}
\label{E:2.4q}
\lim_{t\to\infty}\,\limsup_{N\to\infty}\,
\Bigl|
\frac1t\,\texte^{\sqrt{2\pi}\,t}P\bigl(M_N\ge m_N+t\bigr) -
C_\star \Bigr| = 0 \,.
\end{equation}
Moreover, with $\psi$ as in Theorem~\ref{thm1.7}, for any open set $A\subset[0,1]^2$,
\begin{equation}
\label{E:2.5q}
\lim_{t\to\infty}\,\limsup_{N\to\infty}
\Bigl| P\bigl(N^{-1}\text{\rm argmax}\,h\in A\big|M_N\ge m_N+t\bigr)
- \int_A\psi(x)\textd x \Bigr| = 0 \,.
\end{equation}
\end{theorem}

\begin{remark}
By invoking Corollary~\ref{cor1.5} and the conformal transformation rule \eqref{E:1.20}, we are actually able to conclude that
\begin{equation}
\psi(x)=\frac3\pi\bigl(1-|g(x)|^2\bigr)^2|g'(x)|^{-2},
\end{equation}
where~$g$ is a conformal bijection of $(0,1)^2$ onto the unit disc~$\D$. Since the function $x\mapsto \frac3\pi(1-|x|^2)^2$ is the square of the conformal radius of~$\D$ from~$x$, the latter implies that $\psi$ is, in fact, the square of the conformal radius of~$(0,1)^2$ from~$x$. \textit{Update}: These statements, for general domains, have now all been proved in Biskup and Louidor~\cite{BL2}.
\end{remark}

\subsection{Level sets}
The existing work has not been limited to the maximum $M_N$ only. Indeed, Ding and Zeitouni~\cite{DZ} have also studied the structure of the level sets close to the maximal value (or~$m_N$ above)
\begin{equation}
\label{E:2.7}
\Gamma_N(\lambda):=\{x\in V_N\colon h_x\ge m_N- \lambda\}\,,
\end{equation}
where $\lambda \in \R$. In particular, they derived exponential estimates (in~$\lambda$) on the size of $\Gamma_N(\lambda)$ and controlled distances between the points of~$\Gamma_N(\lambda)$. We restate these results nearly verbatim:

\begin{theorem}[Theorem~1.2 of~\cite{DZ}]
\label{thm-DZ1}
There are constants $0<c<C<\infty$ such that
\begin{equation}
\label{E:2.8a}
\lim_{\lambda\to\infty}\,\liminf_{N\to\infty}\,P\bigl(\texte^{c \lambda}\le|\Gamma_N(\lambda)|\le \texte^{C \lambda}\bigr)
=1.
\end{equation}
\end{theorem}

\begin{theorem}[Theorem~1.1 of~\cite{DZ}]
\label{thm-DZ2}
There is a constant $0<c<\infty$ such that
\begin{equation}
\lim_{r\to\infty}\,\limsup_{N\to\infty}\,P\bigl(\exists u,v\in\Gamma_N(c\log\log r)\colon r\le|u-v|\le N/r\bigr)
=0.
\end{equation}
\end{theorem}
\noindent
The main consequence of these results for us is that they provide the necessary \emph{tightness} for the point processes, both in the $x$-space and $h$-space. This permits us to focus attention on characterizing possible (distributional) limit points.

Broadening our discussion to subjects that are not primarily concerned with the behavior of the extreme points, let us mention also the work of Daviaud~\cite{Daviaud} who studied the size of the set where the field exceeds a constant times~$\log N$. His principal result is that, for $s\in[0,2\sqrt g)$,
\begin{equation}
\label{E:1.5}
\#\bigl\{x\in V_N\colon h_x\ge s\log N\bigr\}=N^{d(s)+o(1)},\quad\text{where}\quad d(s):=2-\frac{s^2}{2g}
\end{equation}
and $o(1)\to0$ in probability. The extreme level sets (at levels of order~$\log N$) thus exhibit a non-trivial fractal structure. The reader should notice the striking similarity with the level sets for~$N^2$ independent Gaussians with variance~$g\log N$.

\smallskip\noindent
\textit{Update}: In a very recent posting, Chatterjee, Dembo and Ding~\cite{Chatterjee-Dembo-Ding} show that an analogous result to \eqref{E:1.5} holds for very general Gaussian fields.
\smallskip

Apart from the Gibbs-Markov property, key to all the above results are the following two technical facts that we record for the purpose of future reference. The first one concerns the asymptotic behavior of the Green function in~$V_N$ away from the boundary while the second one rules out the occurrence of extreme values too close to it. For $\rho \in (0,1)$, we set
\begin{equation}
	V_{N,\rho} := \bigl\{x \in V_N \,\colon \,\, |x-y| > \rho N  \,,\,\, \forall y \in 
	\Z^2\setminus V_N\bigr\} \,.
\end{equation}
Then we have:

\begin{lemma}[Lemma 2.1 of \cite{Daviaud}]
\label{lemma-var-bound}
There exists a universal constant $C > 0$, such that for all $N \geq 1$ and all $x,y \in V_N$,
\begin{equation}
	G_N(x,y) \, \le \, g(\log N - \log (|y-x| \vee 1)) + C \,.
\end{equation}
Moreover for each $\rho \in (0,1) $, there is $C = C(\rho) > 0$ such that for all $N \geq 1$ and $x,y \in V_{N, \rho}$,
\begin{equation}
	\bigl| G_N(x,y) - g(\log N - \log (|y-x| \vee 1)) \bigr| < C \,.
\end{equation}
\end{lemma}

\begin{lemma}[Lemma 3.8 of \cite{BDingZ}]
\label{lemma-extremes-in-bulk}
For all $\lambda > 0$, 
\begin{equation}
	\lim_{\rho \downarrow 0} \limsup_{N \to \infty} \,
		P \bigl(\Gamma_N(\lambda) \setminus V_{N, \rho} \neq \emptyset \bigr) \, = \, 0 \,.
\end{equation}
\end{lemma}

\subsection{Multiplicative chaos and log correlated potentials}
\label{sec2.4}\noindent
Closely related to level sets is the work by Carpentier and Le Doussal~\cite{Carpentier-LeDoussal} (see also Castillo and Le Doussal~\cite{Castillo-LeDoussal}) in the physics literature that (among other things) concerns the $\beta$-dependence of the Gibbs measure
\begin{equation}
\label{E:2.13}
\nu_N(\{x\}):=C_N\,\texte^{\beta h_x},\qquad x\in V_N,
\end{equation}
where~$C_N$ is a normalization constant. The main prediction of \cite{Carpentier-LeDoussal} concerning $\nu_N$ is that the character of this measure changes as $\beta$ increases through a critical value $\beta_\cc:=2/\sqrt{g}$. This is seen rather easily via \eqref{E:1.5} which implies that $\nu_N$ is supported on the level sets from \eqref{E:1.5} with $s\approx s(\beta)$, where $s(\beta):=\beta g$ for $\beta<\beta_\cc$ while $s(\beta):=\beta_\cc g=2\sqrt g$ for $\beta\ge\beta_\cc$. In particular, for~$\beta>\beta_\cc$, $\nu_N$ is supported on the extreme points. 

As noted in~\cite{Carpentier-LeDoussal}, the phenomenon we just described is supposed to occur in all spatial dimensions provided one replaces the DGFF by logarithmically-correlated Gaussian fields. Arguin and Zindy~\cite{AZ} have recently proved this for a one-dimensional model in this class. In fact, they showed that for $\beta>\beta_\cc$, the Gibbs measure remains atomic in the limit $N\to\infty$ and is asymptotically described by a Poisson-Dirichlet law which, we note, is consistent with the extremal point process being Gumbel distributed. (As $\nu_N$ is normalized, the random shift $\alpha^{-1}\log Z$ factors out from the expression.) The analogy with the Random Energy Model is quite striking.

\smallskip\noindent
\textit{Update}: In a recent posting~\cite{Arguin-Zindy2}, Arguin and Zindy have extended some of their conclusions to the DGFF. Pure atomicity, and the characterization in terms of the Poisson-Dirichlet distribution, of the measure~$\beta>\beta_\cc$ is shown in the forthcoming article~\cite{BL3}.
\smallskip

A version of the Gibbs measure \eqref{E:2.13} appears also in the studies of Gaussian multiplicative chaos by Robert and Vargas~\cite{Robert-Vargas}, Allez, Rhodes and Vargas~\cite{ARV} and Duplantier, Rhodes, Sheffield and Vargas~\cite{DRSV1,DRSV2}; see Rhodes and Vargas~\cite{Rhodes-Vargas} for a recent review of this subject going back to Kahane~\cite{Kahane}. In particular, the papers~\cite{DRSV1,DRSV2} are concerned with the construction of a \emph{derivative martingale} which corresponds to the $N\to\infty$ limit of (unnormalized) measure
\begin{equation}
\label{E:2.13'}
\tilde\nu_N(\{x\}):=\bigl[\beta_\cc\text{Var}(h_x)- h_x\bigr]\texte^{\beta_\cc h_x-\frac12\beta_\cc^2\text{Var}(h_x)},\qquad x\in V_N.
\end{equation}
This is an object closely related to the $\beta$-derivative of $\nu_N$ at $\beta:=\beta_\cc$. On the basis of various conjectural statements, the references~\cite{DRSV1,DRSV2} predict the derivative martingale to appear in the place of our~$Z(\textd x)$. In Theorem~\ref{thm1.7} we thus prove a version of this prediction for the DGFF. (The apparent discrepancy in the factor $\frac12$ in front of the variance in the exponents of \eqref{E:1.12t} and \eqref{E:2.13'} comes from the fact that the measure $\textd z$ in \eqref{E:1.12t} has been scaled by a factor of~$K^2$.)

\smallskip\noindent
\textit{Update}: Since the first version of the present paper was circulated, two new papers have appeared dealing with log-correlated Gaussian fields (over continuum space) in any~$d\ge1$. First, a paper of Madaule~\cite{Madaule}, where the Laplace transform representation was shown for the law of the maximum, and a paper by Acosta~\cite{Acosta} where tightness of the maximum was shown under somewhat more general conditions than those of~\cite{Madaule}.
\smallskip

\subsection{Branching Brownian motion}
\label{sec2.5w}\noindent
For the current problem, a very close point of contact with the literature is the Branching Brownian motion (BBM). This is a stochastic process on collections of particles $\{X_i(t)\colon i\in I(t)\}$ that perform independent Brownian motions and undergo splits into two particles at independent exponential waiting times. There is one particle at~$x=0$ at time zero.

The problem originated in the context of nonlinear PDEs with the Kolmogorov, Petrovsky and Piscounov~\cite{KPP} (KPP) --- a.k.a.\ Fisher-Kolmogorov or Fisher-KPP --- equation
\begin{equation}
u_t =\frac12u_{xx} + u^2 - u
\end{equation}
for a function $u=u(t,x)$ where $x\in\R$ and~$t\ge0$. McKean~\cite{McKean} observed that the solution for initial data $u(0, x):=\1_{x\ge0}$ is the probability distribution,
$u(t,x):=P(\,\max_{i\in I(t)} X_i(t)\le x)$, of the particle in BBM that is farthest to the right. Through the works of Kolmogorov, Petrovsky and Piscounov~\cite{KPP} and Bramson~\cite{Bramson1,Bramson2}, it was then established that for
\begin{equation}
m(t):=\sqrt t-\frac3{2\sqrt 2}\log t
\end{equation}
one has $u(x+m(t),t)\to w(x)$ as $t\to\infty$, where $w$ solves the ODE $\frac12w''+\sqrt2 w'+w^2-w=0$ with boundary ``values'' $1$ at $+\infty$ and $0$ at~$-\infty$. Hence, $w$ is the asymptotic distribution function of $\max_{i\in I(t)} X_i(t)-m(t)$ in the limit~$t\to\infty$. Lalley and Sellke~\cite{LS} then proved that $w$ admits a representation by means of the Laplace transform
\begin{equation}
w(x)=E\bigl(\texte^{-c Z\texte^{-\sqrt2\,x}}\bigr),
\end{equation}
where~$c>0$ is a constant and $Z$ is the $t\to\infty$ limit of the so called derivative martingale
\begin{equation}
Z(t):=\sum_{i\in I(t)}\bigl(\sqrt 2\, t- X_i(t)\bigr)\texte^{\sqrt2\, X_i(t)-2t}.
\end{equation}
Notice the remarkable analogy with the statement in \eqref{E:1.9b}.

Building further upon this beautiful structure, Arguin, Bovier and Kistler~\cite{ABK1,ABK2,ABK3} and independently A\"id\'ekon, Berestycki, Brunet and Shi~\cite{ABBS} have recently managed to control the full distribution of the extreme points of the set $\{X_i(t)-m(t)\colon i\in I(t)\}$ as~$t\to\infty$. The limit point process is a cluster process associated with Gumbel law of intensity $\texte^{-\sqrt2\, x}\textd x$, again quite analogously to what we show (modulo the clusters) for the DGFF in Theorem~\ref{thm-main}. A key fact (proved in~\cite{ABK1}) is the separation of time scales: If $X_i(t)$ and~$X_j(t)$ are close to~$m(t)$, then the corresponding Brownian paths split either right at the beginning (i.e., at a time $O(1)$) or or right at the very end (i.e., at a time $t-O(1)$). The splittings in time $O(1)$ give rise to a ``residual'' randomness in the problem; this is the origin of the random variable~$Z$. 

\subsection{Branching Random Walk}
The Branching Random Walk is a discrete-time counterpart of the BBM. It is a process on collections of particles which at integer times (and independently of one another) split into a sample from a random point process --- the \emph{splitting process} --- translated to their position. One is again interested in the behavior of the farthest-to-the-right particle at time~$n$, to be denoted~$M_n$. 

In general, $M_n$ will grow linearly with~$n$ but for properly centered and normalized splitting processes --- the so called boundary cases ---  the growth of~$M_n$ is sublinear. Here McDiarmid~\cite{McDiarmid} showed that the correct order is $\log n$ while Bachmann~\cite{Bachmann} and Bramson and Zeitouni~\cite{BrZ} proved tightness of $M_n-\text{med}(M_n)$ under regularity conditions on the tail of the splitting process. Hu and Shi~\cite{HuShi} and Addario-Berry and Reed~\cite{ABR} then established~$\{M_n-\frac32\log n\colon n\ge1\}$ is tight. Pursuing the strategy that proved to be so useful for the BBM, Biggins and Kyprianou~\cite{BK} showed the convergence of the corresponding derivative martingale~$Z$ while A\"id\'ekon~\cite{Aidekon} established a representation for the limiting law of $M_n-\frac32\log n$ as the Laplace transform of~$Z$.

The Branching Random Walk has played a very important role, both technically and conceptually, in the analysis of the extreme points of the~DGFF. This is thanks to the Gibbs-Markov property of the DGFF.

\section{Main steps of the proof}
\label{sec-3}\noindent
We are now ready to commence the exposition of the proofs. This will be done in three steps which are formulated as separate theorems below. We prove one of these theorems immediately while deferring the proof of the other two to later sections. As we will frequently ``test'' the point processes by integrating them against non-negative functions, let us  write $\langle\eta,f\rangle$ for the integral of~$f$ with respect to~$\eta$.

\subsection{Distributional invariance}
\label{sec2.1}\noindent
Recall that Theorem~\ref{thm-DZ1} implies tightness of the processes $\{\eta_{N,r_N}\}_{N\ge1}$ whenever $r_N\to\infty$ with $r_N/N\to0$. We may thus extract  a subsequence that converges weakly to a process~$\eta$ and try to characterize the limiting distribution.
A key observation is that (any such)~$\eta$ is invariant under ``Dysonization'' of its points by a simple diffusion.

Let $W_t$ denote the standard Brownian motion and recall that $\alpha:=2/\sqrt g$. Given a measurable function $f\colon[0,1]^2 \times \R \to[0,\infty)$, let
\begin{equation}
\label{E:3.1}
f_t(x,h)=-\log E^0\bigl(\texte^{-f(x,h+W_t-\frac\alpha2 t)}\bigr),\qquad t\ge0,
\end{equation}
where~$E^0$ is the expectation with respect to~$W_t$. Then we have:

\begin{theorem}
\label{thm-1}
Let~$\eta$ be any subsequential distributional limit of the processes $\{\eta_{N,r_N}\}_{N\ge1}$, for some $r_N\to\infty$ with $r_N/N\to0$. Then for any continuous $f\colon[0,1]^2\times\R\to[0,\infty)$ with compact support and all $t \geq 0$,
\begin{equation}
\label{E:3.2}
E\bigl(\texte^{-\langle\eta,f\rangle}\bigr)=E\bigl(\texte^{-\langle\eta,f_t\rangle}\bigr)\,.
\end{equation}
\end{theorem}

Let us give some intuition into what \eqref{E:3.2} means probabilistically. Pick a sample~$\eta$ of the limit process. The tightness of the processes implies $\eta(C)<\infty$ a.s.\ for any compact set~$C$, while the tightness of the maxima implies that the support of $\eta$ is bounded from above in the height coordinate. It is also known (Theorem~\ref{thm-DZ1}) that the total mass of $\eta$ is infinite a.s. This permits us to write
$\eta=\sum_{i\in\N}\delta_{(x_i,\,h_i)}$, where $(x_i,h_i)$ enumerates the points in the sample. 

Let $\{W_t^{(i)}\colon i\in\N\}$ denote a collection of independent standard Brownian motions which are independent of $\eta$. Set
\begin{equation}
\label{E:3.3qw}
\eta_t:=\sum_{i\in\N}\delta_{(x_i,\,h_i+W_t^{(i)}-\frac\alpha2t)},\qquad t \ge 0.
\end{equation}
Then, by conditioning on $\eta$, we have for all $f \ge 0$, 
\begin{equation}
E \bigl( \texte^{-\langle\eta,f_t\rangle} \bigr) = E \bigl(\texte^{-\langle\eta_t,f\rangle}\bigr) \,.
\end{equation}
Theorem~\ref{thm-1} then implies that 
\begin{equation}
\label{E:3.22}
\eta_t\laweq\eta\,, \qquad  t \ge 0 \,,
\end{equation}
i.e., the law of~$\eta$ is invariant under the above time evolution of its points. In particular, this implies that $\eta_t$ is finite on every compact set a.s. Since this is not \emph{a priori} guaranteed, for~\eqref{E:3.3qw} and \eqref{E:3.22} to be meaningful, a formal treatment of potential blow-ups is necessary. The formulation~\eqref{E:3.2}, 
together with permitting $\langle\eta,f_t \rangle = \infty$, enables us to avoid this issue altogether.

\smallskip
The formal proof of the theorem is lengthy and is thus relegated to Section~\ref{sec4}. To give some understanding, let us show a back-of-the-envelope calculation that harbors the essential ideas.

The starting point is the Gaussian interpolation, $h\laweq \sqrt{1-s}\,h'+\sqrt{s}\,h''$, where $s\in[0,1]$ and $h'$ and~$h''$ are are two independent copies of~$h$. Setting $s:=t/(g\log N)$ for some~$t>0$ that will be fixed throughout, we get
\begin{equation}
\label{E:3.305}
\begin{aligned}
h&\laweq \sqrt{1-\frac t{g\log N}}\,\,\,h'+\sqrt{\frac{t}{g\log N}}\,\,h'' 
\\
&\;=\,\, h'-\frac12\frac t{g\log N}\,h'+\sqrt{\frac{t}{g\log N}}\,\,h''+o(1),
\end{aligned}
\end{equation}
where we used Taylor's theorem to expand the square root and applied that $\max h'=O(\log N)$ with high probability. 

Let now~$x$ be a point well inside~$V_N$ where either $h_{x}\ge m_N-\lambda$ or~$h_{x}'\ge m_N-\lambda$ for some~$\lambda>0$. Recalling that $\Lambda_r(x)$ is an $r$-neighborhood of~$x$, the fact that $h_z''-h_x''=O(1)$ for~$z\in\Lambda_r(x)$ implies
\begin{equation}
h_z\laweq h_z'-\frac12\frac t{g\log N}\,h_z'+\sqrt{\frac{t}{g\log N}}\,\,h_x''+o(1),
\qquad z\in\Lambda_r(x),
\end{equation}
with high probability. As $h'=O(\log N)$ a similar argument shows that, in fact, $h_z-m_N=O(1)$ and $h_z'-m_N=O(1)$ for $z\in\Lambda_r(x)$. Replacing the second occurrence of $h_z'$ by $m_N+O(1)$ and using that $m_N/(g\log N)=\ffrac\alpha2+o(1)$, we thus get
\begin{equation}
\label{E:3.31}
h_z\laweq h_z'-\frac{\alpha}{2}\,t +\sqrt{\frac{t}{g\log N}}\,\,h_x''+o(1),
\qquad z\in\Lambda_r(x),
\end{equation}
As the local maxima of both~$h$ and~$h'$ in~$\Lambda_r(x)$ are achieved at a unique point, this shows that for~$N$ large, they are in fact achieved at the same point. This puts the local maxima of~$h$ and~$h'$ in one-to-one correspondence.

It remains to properly interpret the expression \eqref{E:3.31}. First note that the term containing~$h_x''$ is  asymptotically distributed as~$W_t$. Next we recall that the local maxima of~$h'$ exceeding~$m_N-\lambda$ are all separated by distance of order~$N$. Therefore, if~$x$ and~$y$ are two such local maxima, $h_x''$ and~$h_y''$ have covariance of order unity. Thanks to the normalization by~$\sqrt{\log N}$, and the independence of~$h'$ and~$h''$, they  can thus be regarded as independent. This means that the extreme local maxima of~$h$ are in distribution close to the local maxima of~$h'$ shifted by independent copies of the random variable~$W_t-\frac\alpha2 t$. This readily translates into~\eqref{E:3.22}. 

\subsection{Extracting a Poisson limit law}
The next step consists in extraction of a Poisson limit law based on the distributional identity provided by Theorem~\ref{thm-1}.  The exact statement is as follows:

\begin{theorem}
\label{thm-2}
Suppose~$\eta$ is a point process on $[0,1]^2\times\R$ such that \eqref{E:3.2}, with~$f_t$ as in \eqref{E:3.1}, holds for some $t>0$ and all continuous $f\colon[0,1]^2\times\R\to[0,\infty)$ with compact support. Assume also that $\eta([0,1]^2\times [0,\infty))<\infty$ and
$\eta([0,1]^2\times \R) > 0$ a.s. Then there is a random Borel measure $Z$ on~$[0,1]^2$, satisfying $Z([0,1]^2) \in (0, \infty)$ a.s., such that
\begin{equation}
\label{E:3.9b}
\eta\laweq\text{\rm PPP}\bigl(Z(\textd x)\otimes\texte^{-\alpha h}\textd h\bigr).
\end{equation}
\end{theorem}

\begin{proofsect}{Proof}
The proof follows very closely arguments from Liggett~\cite{Liggett} who studied invariant laws for point measures that evolve by independent Markov chains. 

We start by adapting the proof of Theorem~4.6 of Liggett~\cite{Liggett} to show that~$\eta$ is a Cox process. Consider the following transition kernel on~$\Omega:=[0,1]^2\times\R$:
\begin{equation}
\cmss P_t\bigl((x,h),A\bigr):=P^0\bigl((x,\,h+W_t-\tfrac\alpha2 t)\in A\bigr).
\end{equation}
Since the probability density of $W_t-\tfrac\alpha2 t$ tends to zero uniformly on~$\R$ as~$t\to\infty$, this kernel has (what we would call) the \emph{uniform dispersivity} property,
\begin{equation}
C\subset\Omega\,\,\,\text{compact}\quad\Rightarrow\quad\sup_{(x,h)\in\Omega}\,\cmss P_t\bigl((x,h),C\bigr)\,\underset{t\to\infty}\longrightarrow\,0.
\end{equation}
In particular, if $f\colon[0,1]^2\times\R\to[0,\infty)$ is continuous with compact support and~$f_t$ is as in \eqref{E:3.1}, then $f_t\to0$ uniformly as~$t\to\infty$. Expanding the log in \eqref{E:3.1}, we thus get that
\begin{equation}
\label{E:3.12}
f_t(x,h)\sim E^0\bigl(1-\texte^{-f(x,h+W_t-\frac\alpha2 t)}\bigr),\qquad t\to\infty,
\end{equation}
uniformly in~$(x,h)$.

Inserting \eqref{E:3.12} into \eqref{E:3.2} and using the uniformity in~\eqref{E:3.12} together with the Bounded Convergence Theorem, we then get
\begin{equation}
\begin{aligned}
E\bigl(\texte^{-\langle\eta,f\rangle}\bigr)&=\lim_{t\to\infty}E\Bigl(\,\exp\Bigl\{-\int E^0\bigl(1-\texte^{-f(x,h+W_t-\frac\alpha2 t)}\bigr)\eta(\textd x\,\textd h)\Bigr\}\Bigr)
\\
&=\lim_{t\to\infty}E\Bigl(\,\exp\Bigl\{-\int \bigl(1-\texte^{-f(x,h)}\bigr)M_t(\textd x\,\textd h)\Bigr\}\Bigr),
\end{aligned}
\end{equation}
where $M_t$ is the random measure on~$\Omega$ that is defined by
\begin{equation}
M_t(A):=\int \eta(\textd x\,\textd h)\,\cmss P_t\bigl((x,h),A\bigr).
\end{equation}
Multiplying~$f$ by~$\lambda>0$ and taking~$\lambda\downarrow0$ for various~$f$'s shows that the family $\{M_t\colon t>0\}$ is tight. We may thus find a sequence $t_n\to\infty$ so that $M_{t_n}\Lawarrow M$, where $M$ is a locally-finite random Borel measure on $[0,1]^2 \times \R$. Hence,
\begin{equation}
E\bigl(\texte^{-\langle\eta,f\rangle}\bigr)
=E\Bigl(\,\exp\Bigl\{-\int \bigl(1-\texte^{-f(x,h)}\bigr)M(\textd x\,\textd h)\Bigr\}\Bigr),
\end{equation}
i.e., $\eta$ has the law of $\text{PPP}(M)$ for some random measure~$M$ on~$\Omega$. Note that
\begin{equation}
\label{E:3.20}
	M([0,1]^2\times [0,\infty))<\infty \ \text{ a.s.} \quad \text{and} \quad
	M([0,1]^2\times \R)>0 \ \text{ a.s.}
\end{equation}
thanks to the Poisson structure and the assumptions on $\eta([0,1]^2 \times [0,\infty))$ 
and $\eta([0,1]^2 \times \R)$.

Using \eqref{E:3.2} one more time now shows $M\cmss P_t\laweq M$ for every~$t>0$. Our aim is to apply Corollary~3.8 of Liggett~\cite{Liggett} which states that if $\cmss P$ is a transition kernel associated with a random walk on a non-compact Abelian group and the kernel has no proper closed invariant subgroups, then 
\begin{equation}
M\cmss P\laweq M\quad\Rightarrow\quad M\cmss P= M\text{ a.s.}
\end{equation}
In our case the requirement of absence of proper closed invariant subgroups fails, because the spatial coordinate is not moved by~$\cmss P_t$. We argue as follows: Let~$A$ be a Borel set in $[0,1]^2$ and define a random measure~$L_A$ on~$\R$ by $L_A(B):=M(A\times B)$. Define a Markov kernel $\cmss Q_t$ on~$\R$ by
\begin{equation}
\cmss Q_t(h,B)=P^0\bigl(h+W_t-\tfrac\alpha2\,t\in B\bigr).
\end{equation}
Then $M\cmss P_t\laweq M$ implies $L_A\cmss Q_t\laweq L_A$. The kernel~$\cmss Q_t$ does satisfy the conditions of Corollary~3.8 of Liggett~\cite{Liggett}, and so $L_A\cmss Q_t=L_A$ a.s.\ for every~$t>0$. Thanks to $M([0,1]^2\times [0,\infty))<\infty$ we also have $L_A([0,\infty))<\infty$ a.s. For the ultimate conclusion we now invoke:

\begin{lemma}
\label{lemma-invar}
Suppose~$\nu$ is a Borel measure on~$\R$ such that $\nu\cmss Q_t=\nu$ for some~$t>0$ and also $\nu([0,\infty))<\infty$. Then $\nu(\textd h)=\alpha\nu([0,\infty))\,\texte^{-\alpha h}\textd h$.
\end{lemma}

\noindent
Postponing the proof of the lemma, we set
\begin{equation}
Z(A):=\alpha L_A\bigl([0,\infty)\bigr)=\alpha M\bigl(A\times[0,\infty)\bigr)
\end{equation}
and use the lemma to infer that $L_A(\textd h)=Z(A)\texte^{-\alpha h}\textd h$ a.s. Hence
\begin{equation}
\label{E:3.18}
M(A\times B)=Z(A)\int_B\texte^{-\alpha h}\,\textd h,\qquad \text{a.s.}(A)
\end{equation}
Although the exceptional set in \eqref{E:3.18} may depend on~$A$, the fact that a countable collection of sets of the form $A\times B$ generates all Borel sets in $[0,1]^2\times\R$ shows that \eqref{E:3.18} defines a Borel measure~$Z(\textd x)$ which obeys $M(\textd x\,\textd h)=Z(\textd x)\texte^{-\alpha h}\textd h$ on a set of probability one.
The fact that $Z([0,1]^2) \in (0, \infty)$ a.s. is a consequence of~\eqref{E:3.20}.
\end{proofsect}

For the proof to be complete, we still need to prove the lemma:
\begin{proofsect}{Proof of Lemma~\ref{lemma-invar}}
Since $\nu\cmss Q_t = \nu * \cmss Q_t(0, \cdot)$ and the support of $\cmss Q_t(0, \cdot)$ is all of~$\R$, we may appeal to Theorem 3' in Deny~\cite{Deny} (based on Choquet-Deny~\cite{ChoquetDeny}). This gives that all solutions to $\nu * \cmss Q_t(0, \cdot)  = \nu$ are  linear combinations with non-negative coefficients of the exponential measures $\nu_\kappa(\textd h):=\texte^{-\kappa h}\textd h$ for $\kappa$ that obey $\int \texte^{\kappa h}  \cmss Q_t(0, \textd h) = 1$. The latter equation permits only $\kappa=0, \alpha$, regardless of the value of $t > 0$. Since $\nu([0,\infty))<\infty$, the coefficient of $\nu_0$ must be $0$ and therefore that of $\nu_\alpha$ must be $\alpha \nu([0, \infty))$.
\end{proofsect}

\subsection{Uniqueness of intensity measure}
\label{sec3.3}\noindent
At this point we know --- and this is, in a sense, the principal contribution of the present paper --- that any subsequential weak limit~$\eta$ of $\{\eta_{N,r_N}\}$ is a Poisson point process with intensity $Z(\textd x)\texte^{-\alpha h}\textd h$, for some random~$Z(\textd x)$ on~$[0,1]^2$.
Our next task is to prove that all such subsequential limits are the same; i.e., that the law of the random measure~$Z(\textd x)$ is unique. In this part we will draw heavily on the techniques developed in Bramson, Ding and Zeitouni~\cite{BDingZ}.

We begin by noting that, by \eqref{E:1.9b}, whenever the maximum converges in law (along a subsequence), the limit of the distribution functions is the Laplace transform of the total mass $Z:=Z([0,1]^2)$ of~$Z(\textd x)$. So, if the maximum of the DGFF centered by~$m_N$ has a unique limit (which we know thanks to the main result of~\cite{BDingZ}), the law of $Z([0,1]^2)$ is unique. In light of this, it will thus come as no surprise that to show the uniqueness of the law of the full measure we need the following generalization of Theorem~\ref{thm2.1}:

\begin{theorem}
\label{thm-3}
Let $(A_1, \dots, A_m)$ be a collection of disjoint  non-empty open subsets of $[0,1]^2$. Then the law of
$\bigl( \max \bigl\{ h_x \colon  x \in V_N,\, \tfrac{x}{N} \in A_l \bigr\} - m_N \bigr)_{l=1}^m$
converges weakly as $N \to \infty$.
\end{theorem}

The proof is a modification of the proof of Theorem~1.1 of~\cite{BDingZ} (which we restated as Theorem~\ref{thm2.1}). However, as the technical details are somewhat demanding, we relegate it to Section~\ref{sec5}. Assuming Theorems~\ref{thm-1} and~\ref{thm-3}, we can complete the proof of our main result:

\begin{proofsect}{Proof of Theorem~\ref{thm-main} from Theorems~\ref{thm-1} and~\ref{thm-3}}
For some~$r_N$ with $r_N\to\infty$ and $r_N/N \to 0$ as $N \to \infty$, let $\eta$ be a sub-sequential  limit of $\{{\eta}_{N, r_N}\}$. Let~$Z(\textd x)$ be the random measure so that \eqref{E:3.9b} holds and let $h^\star_A$ denote the maximum of $h_x$ over $x\in V_N$ with $x/N\in A$. Note that in light of Theorem~\ref{thm-3}, $h_A^\star - m_N$ is tight. Given any collection $(A_1, \dots, A_m)$ of disjoint non-empty open subsets of $[0,1]^2$ such that $Z(\partial A_l)=0$ a.s.\ for each $l=1,\dots,m$, we then have
\begin{equation}
\label{E:3.27}
E \biggl(\exp \Bigl\{ -{\alpha}^{-1} \sum_{l=1}^m \texte^{-\alpha t_l} Z(A_l) \Bigr\}\biggr) = 
\lim_{N \to \infty} P \bigl(\,
	 h^\star_{A_l} - m_{N} \leq t_l\,,\,\, l=1, \dots, m \bigr)
\end{equation}
for any $t_1, \dots, t_m \in \R$. (The convergence in Theorem~\ref{thm-3} ensures this for a dense set of $t_l$'s; the continuity of the left-hand side then extends this to all of~$\R^m$.)

By Theorem~\ref{thm-3} again, the right-hand side of \eqref{E:3.27} is the same for all subsequences and so this proves uniqueness of the law of~integrals with respect to~$Z(\textd x)$ of all positive simple functions on open sets $A$ with $Z(\partial A)=0$ a.s. Using standard arguments this can be extended to the class of all continuous functions on $[0,1]^2$. Hence, the law of $Z(\textd x)$ is itself unique.

The a.s.\ finiteness of the total mass of~$Z$ arises (via the arguments in Theorem~\ref{thm-2}) from the tightness of the upper tail of~$\eta$. The fact that $Z(A)>0$ a.s.\ for any $A\subset[0,1]^2$ open is a consequence of the fact that~$\eta(A\times\R)>0$ a.s.\ for all such~$A$. This follows from the Gibbs-Markov property of the DGFF and the fact that~$A$ contains an open square.
\end{proofsect}

Similar arguments as used above also permit us to establish an explicit form of the limiting law of the maximum, along with its scaled position:

\begin{proofsect}{Proof of Corollary~\ref{cor1.5}}
Let~$A$ be an open set with $Z(\partial A)=0$ a.s.\ and recall the notation $h^\star_A$ from the previous proof. Note that \eqref{E:3.27} and a continuity argument imply that the joint law of $(h^\star_A-m_N,\,h^\star_{A^\cc}-m_N)$ converges to the corresponding object for the limit point process. Since
\begin{equation}
\nu_N\bigl(A\times(-\infty,t]\bigr)=P\bigl(\,h^\star_A - m_{N}\le t,\, h^\star_{A^\cc}-m_N< h^\star_A-m_N\Bigr),
\end{equation}
the claim follows by a straightforward calculation.
\end{proofsect}

With Theorem~\ref{thm-main} in hand, we are also able to establish the basic properties of the $Z$-measure as stated in Theorem~\ref{thm1.3}. First we prove a lemma:

\begin{lemma}
\label{lemma3.6}
Let~$Z:=Z([0,1]^2)$ be related to a (subsequential) limit~$\eta$ of $\{\eta_{N,r_N}\}_{N\ge1}$ as in Theorem~\ref{thm-2}. Let~$C_\star$ denote the limit value in \eqref{E:2.4q}. Then
\begin{equation}
\label{E:3.29}
\lim_{\lambda\downarrow0}\,\frac{E\bigl(\,1-\texte^{-\lambda Z}\bigr)}{\lambda\log(\ffrac1\lambda)}= C_\star
\end{equation}
and
\begin{equation}
\label{E:3.30}
\lim_{\lambda\downarrow0}\,\frac{E(\,Z\texte^{-\lambda Z})}{\log(\ffrac1\lambda)}= C_\star.
\end{equation}
\end{lemma}

\begin{proofsect}{Proof}
Let $t$ and~$\lambda$ be related via $\lambda=\alpha^{-1}\texte^{-\alpha t}$ and note that $\alpha=\sqrt{2\pi}$ implies $\lambda\log(\ffrac1\lambda)=t\texte^{-\sqrt{2\pi}\,t} (1+o(1))$. The Laplace-transform representation \eqref{E:1.9b} of the limit law of the maximum 
and the asymptotic formula \eqref{E:2.4q} then readily yield the asymptotic expression \eqref{E:3.29}.

For \eqref{E:3.30}, we first note note that
\begin{equation}
E\bigl(\,1-\texte^{-\lambda Z}\bigr)=\int_0^\lambda E\bigl(\,Z\texte^{-\lambda' Z}\bigr)\textd\lambda'.
\end{equation}
Integrating instead from $\theta \lambda$ to $\lambda$ where $\theta \in (0,1)$ and using monotonicity we then get
\begin{equation}
\lambda(1-\theta)E\bigl(\,Z\texte^{-\lambda Z}\bigr)\le
E\bigl(\,\texte^{-\lambda\theta Z}-\texte^{-\lambda Z}\bigr)\le\lambda(1-\theta)E\bigl(\,Z\texte^{-\lambda\theta Z}\bigr).
\end{equation}
Dividing by $\lambda\log(\ffrac1\lambda)$ and letting $\lambda \downarrow 0$, the middle expression converges to $(1-\theta) C_\star$ by~\eqref{E:3.29} and the left-hand inequality then gives an upper bound on the \emph{limes superior} in \eqref{E:3.30}. Dividing instead by $\lambda\log(\ffrac1{\theta \lambda}) = \lambda\log(\ffrac1\lambda)(1+o(1))$ and using now the right-hand inequality gives the corresponding lower bound on the \emph{limes inferior} in~\eqref{E:3.30}.
\end{proofsect}

\begin{proofsect}{Proof of Theorem~\ref{thm1.3} from Theorem~\ref{thm-main}}
Throughout we assume that $\eta_{N,r_N}\Lawarrow\eta$ with representation using the $Z(\textd x)$ measure with total mass~$Z$. Equation \eqref{E:3.30} and a standard Tauberian theorem applied to the monotone function $t \mapsto E(Z \1_{\{Z \leq t\}})$ (c.f., e.g., Theorem 2 in section XIII.5 of Feller~\cite{Feller}) imply
\begin{equation}
E(Z\1_{\{Z\le t\}})=[C_\star+o(1)]\log(t),\qquad t\to\infty.
\end{equation}
Using Fubini-Tonelli to write, for $q>0$,
\begin{equation}
E\Bigl(\frac Z{(\log (Z\vee 2))^q}\Bigr)=\int_2^\infty \frac q{t(\log t)^{q+1}} E(Z\1_{\{Z\le t\}})\,\textd t,
\end{equation}
the integrability claims for the positive powers follow. 

For negative powers we first note that \eqref{E:2.3q} yields
\begin{equation}
E(\texte^{-\lambda Z})\le C\texte^{-c\lambda^c}
\end{equation}
for some $0<c<C<\infty$ and $\lambda$ large. This implies $P(Z<t)\le C\texte^{-ct^{-c}}$ and so $Z^{-1}$ has even stretched exponential moments.

Moving to non-atomicity of~$Z(\textd x)$, a key observation is that if~$Z$ had an atom of strength~$c>0$ in $[0,1]^2$, then for each $\epsilon > 0$, one could find~$\lambda\in\R$ such that, conditional on $Z$, the following would hold with probability at least $1-\epsilon$: For each $\delta>0$ there is a ball~$B$ of radius~$\delta>0$ such that $\eta\big((B \cap [0,1]^2) \times[\lambda,\infty) \big) \ge 2$. By Theorem~\ref{thm-DZ2}, this can be ruled out (uniformly in~$N$) for processes~$\eta_{N,r_N}$ but work is needed to extend this to the limit process~$\eta$. Throughout we will regard~$\eta_{N,r_N}$ and~$\eta$ as measures on all of~$\R^2$.

First we claim that for each~$\delta>0$ there is~$n\ge1$ and balls~$B_1,\dots,B_n$ with centers in~$[0,1]^2$ and radii in~$(\delta,2\delta)$ such that
\begin{enumerate}
\item[(1)] every ball~$B$ of radius~$\delta$ centered in~$[0,1]^2$ is entirely contained in one of~$B_1,\dots,B_n$, and
\item[(2)] for each $\lambda\in\R$,
\begin{equation}
\label{E:3.28q}
\bigl(\eta_{N,r_N}(B_i\times [\lambda,\infty))\colon l=i,\dots,m\bigr)\,\underset{N\to\infty}\Lawlongarrow\,\bigl(\eta(B_i\times [\lambda,\infty))\colon i=1,\dots,m\bigr)
\end{equation}
as random vectors on~$\R^m$.
\end{enumerate}
To see this, we observe that \eqref{E:3.28q} holds for any open sets~$B_1,\dots,B_n$ as long as they obey $P(Z(\partial B_i)=0)=1$ for all~$i=1,\dots,n$. (This is thanks to $\eta_{N,r_N}\Lawarrow\eta$ and the continuity of the intensity measure in~$h$-variable.) So one thus first chooses a cover $B_1',\dots,B_n'$ of~$[-\delta,1+\delta]^2$ by balls of radii in $\frac32\delta$ so that~(1) holds and then increases the radius of each ball slightly so that the $Z$-measure of their boundary vanishes a.s.

Now fix~$\epsilon>0$ and~$\lambda\in\R$. By Theorem~\ref{thm-DZ2}, there is $\delta>0$ such that the probability of there existing a ball $B$ of radius $2\delta$ with $\eta_{N,r_N}\big( (B \cap [0,1]^2) \times[\lambda,\infty) \big)\ge2$ is less than~$\epsilon$ once~$N$ is sufficiently large. Using this~$\delta$ for the above cover yields
\begin{equation}
\label{E:3.29a}
P\bigl(\max_{1\le l\le m}\eta(B_l\times [\lambda,\infty))\ge2\bigr)\le\epsilon.
\end{equation}
However, if~$F_c$ is the event that~$Z(\textd x)$ has an atom of strength at least~$c>0$, by conditioning on~$Z(\textd x)$ and a simple calculation for the Poisson law we get
\begin{equation}
\label{E:3.32xx}
P\bigl(\,\max_{1\le l\le m}\eta(B_l\times [\lambda,\infty))\le1\,\big|\, F_c\bigr)\le
[1+c\alpha^{-1}\texte^{-\alpha\lambda}]\exp\{-c\alpha^{-1}\texte^{-\alpha\lambda}\},
\end{equation}
Denoting the the right-hand side of \eqref{E:3.32xx} by~$\theta$, from \eqref{E:3.29a} we infer $P(F_c)\le\epsilon/(1-\theta)$. This can be made arbitrarily small by adjusting~$\epsilon$ and~$\lambda$ (which have insofar been arbitrary) and so the measure~$Z(\textd x)$ is non-atomic after all.

A completely analogous argument, albeit based on Lemma~\ref{lemma-extremes-in-bulk} instead of Theorem~\ref{thm-DZ2}, proves that $Z(\textd x)$ puts zero mass on the boundary of the unit box.
\end{proofsect}

At this point we have established Theorems~\ref{thm-main} and~\ref{thm1.3} and Corollary~\ref{cor1.5} subject to  Theorems~\ref{thm-1} and~\ref{thm-3}. We will now move on to proving these theorems as well.
Theorem~\ref{thm1.7} requires the structure underlying the proof of Theorem~\ref{thm-3} and so we defer its proof to Section~\ref{sec5.1}.

\section{Distributional invariance}
\label{sec4}\noindent
The first item left to prove from the previous section is the distributional identity from Theorem~\ref{thm-1} that should hold for any subsequential limit of the extremal point processes. 

\subsection{Tightness revisited}
\label{sec4.1}\noindent
As remarked earlier, the statements in Theorem~\ref{thm-DZ1} provide tightness for the corresponding extreme-value point process. However, our proof of Theorem~\ref{thm-1} requires a bound that is summable in $\lambda$ (uniformly in~$N$) instead of just the convergence to one in~\eqref{E:2.8a}. 

\begin{proposition}
\label{prop-tight}
There exist $\beta > 0$ such that for all~$\lambda>1$ and all large-enough $\kappa>0$, 
\begin{equation}  
\label{eqn:17}
\sup_{N\ge1}\,\,P\bigl(|\Gamma_N(\lambda)| > \texte^{\kappa \lambda}\bigr) \le \texte^{-\beta \kappa \,\lambda} \,.
\end{equation}
In particular, every weak subsequential limit~$\eta$ of the processes $\{\eta_{N,r_N}\}$
(for any choice of $r_N$) satisfies,
\begin{equation}
P\bigl(\eta([0,1]^2\times[-\lambda,\infty)) > \texte^{\kappa \lambda}\bigr)\le\texte^{-\beta\kappa \lambda}
\end{equation}
for all~$\lambda>1$ and all~$\kappa$ sufficiently large.
\end{proposition}

The proof of the proposition will be based on two lemmas. Recall the following object from Ding and Zeitouni~\cite{DZ}:
\begin{equation}
\label{E:4.3}
S_{\ell,N} =  S_{\ell,N}(h) :=\max\Bigl\{\sum_{x\in U}h_x\colon U\subset V_N,\,|U|=\ell\Bigr\} 
\end{equation}
and let $U_{\ell, N} = U_{\ell,N}(h)$ denote the set~$U$ that achieves the maximum. 
Then we have:

\begin{lemma}
\label{lemmaA}
There are constants $c_1,c_2\in(0,\infty)$ such that for all $\lambda>0$,~$\ell\ge1$ and~$N\ge1$,
\begin{equation}
P\bigl(S_{\ell,N}<\ell (m_N-\lambda)\bigr)>c_1-c_2\,\frac \lambda{\log\ell}.
\end{equation}
\end{lemma}

\begin{lemma}
\label{lemmaB}
There are constants $c_3,c_4\in(0,\infty)$ such that the following is true for all $\epsilon,\delta>0$ that are sufficiently small: Set $\sigma:=\epsilon^{-1}\log(1/\delta)$ and, given $\lambda>0$ and $\ell\in\N$, define 
\begin{equation}
\label{E:4.5}
\lambda':=\lambda+c_4\bigl(1+\log \sigma+\sqrt{\log(\ell/\delta)\log \sigma}\bigr).
\end{equation}
If for some $N'\ge1$
\begin{equation}
\label{E:4.4}
P\bigl(S_{\ell,N'}<\ell (m_{N'}-\lambda')\bigr)>\delta,
\end{equation}
then 
\begin{equation}
P\bigl(|\Gamma_{N}(\lambda)|>\ell\bigr)<\epsilon
\end{equation}
holds for all $N\in\N$ with
\begin{equation}
\label{E:4.6b}
N\le c_3 \sigma^{-1/2}\,N'
\,.
\end{equation}
\end{lemma}

Before giving a proof of these lemmas, let us see how they imply Proposition~\ref{prop-tight}.

\begin{proofsect}{Proof of Proposition~\ref{prop-tight} from Lemmas~\ref{lemmaA} and~\ref{lemmaB}}
Let $c_1,c_2$ be the constants from Lemma~\ref{lemmaA}. Given~$\lambda>1$ and $\kappa>1$, define $\ell:=\lceil\texte^{\kappa \lambda}\rceil$ and set $\epsilon:=\texte^{-\beta\kappa\lambda}$ for some~$\beta\in(0,1)$ to be determined. Fix~$\delta\in(0,c_1)$ and let $\lambda'$ be defined as in \eqref{E:4.5}. A simple estimate then shows that for $\beta > 0$ small enough, any $\lambda > 0$ and $\kappa > \kappa_0(\beta, \lambda)$ --- to make $\beta\kappa\lambda$ large --- we have
\begin{equation}
\lambda'\le\lambda\bigl[1+\tilde c_1\sqrt\beta\,\kappa\bigr].
\end{equation}
But then $\lambda'/\log\ell=O(\sqrt\beta)+O(\kappa^{-1})$ and so for~$\beta$ small and~$\kappa$ large we have
\begin{equation}
c_1-c_2\frac{\lambda'}{\log\ell}>\delta.
\end{equation}
By Lemma~\ref{lemmaA}, we get \eqref{E:4.4} for all~$N'\ge1$. Lemma~\ref{lemmaB} then yields
\begin{equation}
\label{E:4.7}
P\bigl(|\Gamma_{N}(\lambda)|>\texte^{\kappa\lambda}\bigr)<\texte^{-\beta\kappa\lambda}
\end{equation} 
for all~$N$ satisfying \eqref{E:4.6b}. But~$N'$ can be taken arbitrary large and so the inequality holds for all~$N\ge1$ as claimed.
\end{proofsect}

We will now move to the proofs of the lemmas. We will plug into an argument from Ding and Zeitouni~\cite{DZ}. Specifically, we will need the following facts:
\settowidth{\leftmargini}{(1111)}
\begin{enumerate}
\item[(1)] Formula (58) of~\cite{DZ} gives
\begin{equation}
\label{E:4.9}
E S_{\ell,N}\le \ell(m_N-c\log\ell)
\end{equation}
for some constant~$c>0$.
\item[(2)] Let $\Xi_{N,r}:=\{(u,v)\in V_N\times V_N\colon r\le|u-v|\le N/r\}$. Lemma~4.6 of~\cite{DZ} states
\begin{multline}
\label{E:4.10}
\qquad\quad
P\bigl(\exists A\subset\Xi_{N,r}\colon|A|\ge\log r,\, (u,v)\in A\,\Rightarrow\,h_u+h_v\ge2[m_N-\lambda\log\log r]\bigr)
\\
\ge1-C\texte^{-\texte^{c\lambda\log\log r}}
\quad
\end{multline}
for some $C,c>0$, all $N\ge1$ and all $r,\lambda\ge C$.
\end{enumerate}
Based on these we can now prove Lemma~\ref{lemmaA}:

\begin{proofsect}{Proof of Lemma~\ref{lemmaA}}
Abbreviate
\begin{equation}
\widetilde S_{\ell,N}:=m_N-\frac1\ell S_{\ell,N}.
\end{equation}
Our goal is to prove a lower bound on the upper tail of~$\widetilde S_{\ell,N}$. For this let us note that on $\{\widetilde S_{\ell,N}>s\}$, the average of the field on any set of size~$\ell$ is less than $m_N-s$. Setting $\ell:=\sqrt{\log r}$, \eqref{E:4.10} shows that with probability at least $1-C\texte^{-2c\lambda\log\ell}$ there is a set of size~$\ell$, where~$h$ averages to at least $m_N-2\lambda\log\ell$. Writing~$s$ in place of~$2\lambda\log\ell$, we thus get
\begin{equation}
\label{E:4.12}
P(\widetilde S_{\ell,N}>s)\le C\texte^{-\texte^{cs}},\qquad \ell\ge C,\, s\ge C'\log\ell.
\end{equation}
The restrictions on $\ell$ and~$s$ on the right come from the requirement in~\eqref{E:4.10} that $r,\lambda \ge C$.

The rest of the proof is a calculation. Consider the expectation $E\widetilde S_{\ell,N}^+$ where the plus denotes the positive part. Write this as the integral of $P(\widetilde S_{\ell,N}>s)$ over~$s\ge0$ and divide the integration domain into three parts: $[0,\lambda)$,  $[\lambda,C'\log\ell)$ and $[C'\log\ell,\infty)$.
In the integral over the last part we may bound the integrand as in \eqref{E:4.12}. Thus we get
\begin{equation}
\label{E:4.13}
E\widetilde S_{\ell,N}^+\le \lambda +C'(\log\ell)\,P(\widetilde S_{\ell,N}>\lambda)+\int_{C'\log\ell}^\infty C\texte^{-\texte^{cs}}\textd s.
\end{equation}
As is easy to check, the last integral is bounded by $C\texte^{-\ell^c}$. 

Now, in order to get a lower bound on $P(\widetilde S_{\ell,N}>\lambda)$ from \eqref{E:4.13}, it suffices to get a suitable lower bound on $E\widetilde S_{\ell,N}^+$. Here we note that, by \eqref{E:4.9},
\begin{equation}
\label{E:4.14}
E\widetilde S_{\ell,N}^+\ge E\widetilde S_{\ell,N}\ge c\log\ell.
\end{equation}
Putting \twoeqref{E:4.13}{E:4.14} together we have
\begin{equation}
P(\widetilde S_{\ell,N}>\lambda)\ge \frac{c\log\ell-\lambda-C\texte^{-\ell^c}}{C'\log\ell}.
\end{equation}
For properly chosen $c_1, c_2 > 0$, this has the form in the statement of the lemma.
\end{proofsect}

\begin{proofsect}{Proof of of Lemma~\ref{lemmaB}}
We will prove the counterpositive. Fix~$\epsilon,\delta>0$ small enough and suppose that for some $\lambda>0$, $\ell\ge1$ and~$N\ge1$,
\begin{equation}
\label{E:4.16}
P\bigl(|\Gamma_{N}(\lambda)|>\ell\bigr) \ge \epsilon.
\end{equation}
Notice that this implies
\begin{equation}
\label{E:4.17}
P\bigl(S_{\ell,N}>\ell(m_{N}-\lambda)\bigr)>\epsilon.
\end{equation}
Let us pick $\rho\in(0,1)$ and define $N'':=\lfloor\rho^{-1}N\rfloor$ and $N':=\lfloor\sigma^{1/2}\rho^{-2}N\rfloor$. Our strategy is to use the method introduced by Ding~\cite{Ding} (called there ``sprinkling'') and used in a similar context by Ding and Zeitouni~\cite{DZ}. To this end, we identify inside~$V_{N'}$ altogether $n:=(\lfloor N'/(4N'')\rfloor)^2$ disjoint translates $V_{N''}^i$, $i=1,\dots,n$, of~$V_{N''}$ such that any two of these are at least distance~$2$ apart and distance $N'/4$ away from the the boundary of $V_{N'}$. In each $V_{N''}^i$ we then fix a translate~$V_{N}^i$ of~$V_N$ which is centered, as much as possible, near the center of~$V_{N''}^i$.

As a first step, consider a box $V_{N}$ of side $N$ contained in, and centered at the same point near the center of, a box~$V_{N''}$ of side~$N''$. If $h''$ is a DGFF on $V_{N''}$ with zero boundary conditions, then we claim that
\begin{equation}
\label{E:4.18}
P\bigl(S_{\ell,N}(h'')>\ell(m_{N}- \lambda )\bigr)>\frac\epsilon2 
\end{equation}
(notice that $h''$ is defined on $V_{N''}$ but $S_{l,N}(h'')$ is restricted to subsets of $V_N$).
To see this, note that the Gibbs-Markov property implies that in~$V_{N}$ we have $h''\laweq h+\varphi''$, where~$h$ is the DGGF in~$V_{N}$ and $\varphi''$ is a mean-zero Gaussian field that is independent of~$h$. Therefore, for any $\lambda$, 
\begin{equation}
\bigl\{S_{\ell,N}(h'')>\ell(m_{N}- \lambda)\bigr\}\supseteq \bigl\{S_{\ell,N}(h)>\ell(m_{N}- \lambda)\bigr\}\cap\Bigl\{\sum_{x\in U_{\ell,N}(h)}\varphi''(x) \geq 0\Bigr\},
\end{equation}
where $U_{\ell,N}$ is defined below~\eqref{E:4.3}. The independence of~$\varphi''$ of~$h$ then ensures that, conditional on~$h$, the probability of the last event is~$\ffrac12$. Using \eqref{E:4.17}, the bound \eqref{E:4.18} follows.

We may now decompose the DGFF $h'$ on $V_{N'}$ in the usual way, by conditioning on its value at the boundary of the boxes $V_{N''}^i$, for $i=1, \dots, N$. Under this decomposition, $h'$ restricted to each box $V_{N''}^i$ is the sum $h''_i + \varphi'$, where~$h''_i$ is the DGFF in~$V_{N''}^i$ and $\varphi'$ is the ``binding field'', i.e. $\varphi' =  E \bigr( h' \ |\  h'_x \colon  x \in \partial V_{N''}^i \,,\,\, i=1, \dots, n \bigl)$. Moreover, $\varphi'$ and all $h''_i$ are independent, $\varphi'$ has zero mean and for all $x \in \bigcup_{i=1}^n V_N^i$ 
\begin{equation}
\label{E:4.20}
\text{Var}(\varphi'(x))\le C''\log(N'/N'') \,.
\end{equation}
This is a direct consequence of Lemma~\ref{lemma-var-bound} and the enforced distance of $x$ from the boundary of both $V_{N'}$ and the corresponding box $V_{N''}^i$.

Denoting by $K(h')$ the index~$i$ where $S_{\ell,N}(h''_i)$ is maximized, we then get for all $a,b \in \R$,
\begin{multline}
\label{E:4.21}
\qquad
\bigl\{S_{\ell,N'}>\ell(m_{N'}- a-b)\bigr\}
\\
\supseteq\bigl\{\max_{1\le i\le n}S_{\ell,N}(h''_i)>\ell(m_{N'}- a)\bigr\}
\cap\Bigl\{\,\min_{x\in U_{\ell,N}(h''_{K(h')})}\varphi'(x)\ge-b\Bigr\}.
\qquad
\end{multline}
Set $b :=c_4\sqrt{\log(\ell/\delta)\log\sigma}$, for some $c_4 > 0$ to be defined shortly. By the union bound, a standard Gaussian estimate, the estimate \eqref{E:4.20} and the fact that $\log(N'/N'')$ is order $\log\sigma$, the probability of the last event in~\eqref{E:4.21} is
\begin{equation}
\label{E:4.22}
P\Bigl(\,\min_{x\in U_{\ell,N}(h''_{K(h')})}\varphi'(x)\ge-b \Bigr)
\ge 1-\ell\texte^{-cb^2/\log\sigma}  \,.
\end{equation}
The right hand side above can be lower bounded by $1-\delta/2$ uniformly in $\ell$, if $c_4$ is chosen large enough.

For the first term on the right-hand side of~\eqref{E:4.21} we set $a :=m_{N'}-m_N+\lambda$ and assume that~$\rho$ so small that
\begin{equation}
\label{E:4.23}
P\Bigl(\,\max_{1\le i\le n}S_{\ell,N}(h''_i)>\ell(m_{N}- \lambda)\Bigr)\ge1-(1-\epsilon/2)^n 
\geq 1 - \delta/2 \,.
\end{equation}
Since $m_N-m_{N'} \le c\log\sigma$ for $\delta, \epsilon$ small, it follows that $\lambda'$ as defined in~\eqref{E:4.5} satisfies $\lambda' > a+b$. Combining \twoeqref{E:4.21}{E:4.23} we  conclude 
\begin{equation}
\label{E:4.27a}
P\bigl(S_{\ell,N'}>\ell(m_{N'}-\lambda')\bigr) \ge
P\bigl(S_{\ell,N'}>\ell(m_{N'}-a - b)\bigr)
\ge 1-\delta \,.
\end{equation}
This is the converse of~\eqref{E:4.4} as we aimed to show. 
\end{proofsect}

\subsection{Proof of distributional invariance}
We are now ready to move to the proof of Theorem~\ref{thm-1}. Let~$h',h''$ be independent copies of the DGFF in~$V_N$. Let $t>0$ be a fixed number which does not change with $N$ and abbreviate
\begin{equation}
	\hat{h}' := \sqrt{1-\frac t{g\log N}}\,\,\,h'
	\quad\text{and}\quad 
	\hat{h}'' := \sqrt{\frac{t}{g\log N}}\,\,\,h''.
\end{equation}
As mentioned in Section~\ref{sec2.1}, key to our proof is the fact that the field
\begin{equation}
\label{E:4.28a}
h :=\hat{h}' + \hat{h}''
\end{equation}
also has the law of the DGFF in~$V_N$. In what follows we will regard $h'$ and~$h''$ as realized on the same probability space and~$h$ defined by \eqref{E:4.28a}.

Recall our notation~$\Gamma_N(\lambda)$ from \eqref{E:2.7} for the level set~``$\lambda$ units below~$m_N$'' and let $\Theta_{N,r}$ denote the set of all $r$-local maxima,
\begin{equation}
\Theta_{N,r}:=\bigl\{x\in V_N\colon h_x=\max_{z\in\Lambda_r(x)}h_z\bigr\}.
\end{equation}
Clearly, for any $f\colon[0,1]^2\times\R\to[0,\infty)$,
\begin{equation}
\label{E:4.31}
\langle{\eta}_{N,r},f\rangle=\sum_{x\in\Theta_{N,r}}f\bigl(\tfrac xN,h_x-m_N\bigr).
\end{equation}
If~$f$ is compactly supported, we can insist on $x\in\Gamma_N(\lambda)$ in the sum as well, provided~$\lambda$ is so large that~$f$ vanishes on $[0,1]^2\times(-\infty,-\lambda)$. Our first observation is that there is not much difference in using~$r$ or~$r_N$ in the definition of local maxima:

\begin{lemma}
\label{lemma-T1}
Let $f\colon[0,1]^2\times\R\to[0,\infty)$ be a measurable function with compact support.
For any $r_N$ with $r_N\to\infty$ and $r_N/N\to0$,
\begin{equation}
\lim_{r\to\infty}\limsup_{N\to\infty}\,P\bigl(\langle{\eta}_{N,r_N},f\rangle \neq \langle{\eta}_{N,r},f\rangle\bigr)=0.
\end{equation}
\end{lemma}

\begin{proofsect}{Proof}
Pick $\lambda>0$. We will first prove
\begin{equation}
\label{E:4.31a}
\lim_{r\to\infty}\limsup_{N\to\infty}\,P\Bigl(\Gamma_N(\lambda)\cap\bigl(\Theta_{N,r}\,{\scriptstyle\triangle}\,\Theta_{N,r_N}\bigr)\ne\emptyset\Bigr)=0.
\end{equation}
Indeed, once $r_N\ge r$ we have $\Theta_{N,r_N}\subset\Theta_{N,r}$ and so, on the event in question, $\Gamma_N(\lambda)$ contains two vertices~$u$ and~$v$ with $r\le|u-v|\le r_N$. Since also~$r_N\le N/r$ for~$N$ large, the probability of such an event tends to zero by Theorem~\ref{thm-DZ2} in the limit $N\to\infty$ followed by $r\to\infty$.

Returning to the main statement, we now set~$\lambda$ to be so large that $f$ vanishes on $[0,1]^2\times(-\infty,-\lambda)$. Then $\langle{\eta}_{N,r_N},f\rangle-\langle{\eta}_{N,r},f\rangle\ne0$ implies containment in the event in \eqref{E:4.31a} and so the corresponding probability tends to zero as well.
\end{proofsect}

Lemma~\ref{lemma-T1} permits us to focus on $\langle{\eta}_{N,r},f\rangle$ for a fixed~$r$. Let us use $\Gamma_N'(\lambda)$, resp., $\Theta_{N,r}'$ to denote the same objects as $\Gamma_N(\lambda)$, resp., $\Theta_{N,r}$, except that they are defined for~$h'$ instead of~$h$. A main technical step is the replacement of condition $x\in\Theta_{N,r}$ in the sum defining $\langle{\eta}_{N,r},f\rangle$ by $x\in\Theta'_{N,r}$. We formulate this as follows:

\begin{proposition}
\label{lemma-T3}
Let $f\colon[0,1]^2\times\R\to[0,\infty)$ be continuous with compact support. For any~$\epsilon>0$,
\begin{equation}
\label{E:4.48}
\lim_{r\to\infty}\limsup_{N\to\infty}\,P\biggl(\,\Bigl|\,\sum_{x\in \Theta_{N,r}}
\!\!\!f\bigl(\tfrac xN,h_x-m_N\bigr)\,-\sum_{x\in\Theta_{N,r}'}\!\!\!f\bigl(\tfrac xN,h_x-m_N\bigr)\Bigr|>\epsilon\biggr)=0.
\end{equation}
\end{proposition}

We will prove this by showing that the values of~$f$ on $\Theta_{N,r}$, restricted to a proper level set of~$h$, are controllably close to the values of~$f$ on $\Theta_{N,r}'$, restricted to a proper level set of~$h'$. Since~$f$ is compactly supported, these level sets can be chosen so that all terms effectively contributing to the sums in \eqref{E:4.48} are accounted for.
Our first step is a control of the difference between the level sets of~$h$ and~$h'$: 

\begin{lemma}
\label{lemma-T2}
We have
\begin{equation}
\lim_{\lambda\to\infty}\liminf_{N\to\infty}\,P\bigl(\,\Gamma_N'(\lambda) \subseteq \Gamma_N(2\lambda)\bigr) = 1
\end{equation}
and, similarly,
\begin{equation}
\lim_{\lambda\to\infty}\liminf_{N\to\infty}\,P\bigl(\,\Gamma_N(\lambda) \subseteq \Gamma'_N(2\lambda)\bigr) = 1.
\end{equation}
\end{lemma}

\begin{proofsect}{Proof}
For the first event we note that, if $h_x'\ge m_N-\lambda$ with $\lambda>0$, then $\hat h_x'\ge m_N-\lambda-c$, for some~$c>0$ that depends only on~$t$. Hence, $\hat h_x'+\hat h_x''\le m_N-2\lambda$ implies $\hat h_x''\le c-\lambda$ and so
\begin{equation}
\Gamma_N'(\lambda)\setminus\Gamma_N(2\lambda)
\subseteq\Bigl\{x\in V_N\colon h_x'\ge m_N-\lambda,\,\hat h_x'' <  c-\lambda\Bigr\}.
\end{equation}
Since $h''$ is independent of~$h'$, has zero mean and --- as implied by Lemma~\ref{lemma-var-bound} --- satisfies $\Var(\hat h''_x) < C$, once $\lambda>c$ we get
\begin{equation}
P\Bigl(\,\Gamma_N'(\lambda)\setminus\Gamma_N(2\lambda)\ne\emptyset\Bigr)
\le P\bigl(|\Gamma_N'(\lambda)|>\texte^{C'\lambda}\bigr)+\texte^{C'\lambda}\texte^{-(\lambda-c)^2/(2C)}
\end{equation}
by conditioning on the size of $\Gamma_N'(\lambda)$ and a union bound. Theorem~\ref{thm-DZ2} now tells us that, for~$C$ large enough, the first term tends to zero as $N\to\infty$ and~$\lambda\to\infty$. Hence, the left-hand side also vanishes in the stated limit.

For the second event we abbreviate $a_N:=(1-t/(g\log N))^{1/2}$ and notice that, for $n\ge2$ and $x\in V_N$ such that $m_N-(n+1)\lambda<h_x'\le m_N-n\lambda$ we have $\hat h_x'\le m_N-a_Nn\lambda$. Under such conditions, $\hat h_x'+\hat h_x''\ge m_N-\lambda$ implies $\hat h_x''\ge(a_Nn-1)\lambda$. It follows that
\begin{equation}
\label{E:4.35}
\Gamma_N(\lambda)\setminus\Gamma_N'(2\lambda)\subseteq\bigcup_{n\ge2}
\Bigl\{x\in V_N\colon h_x' >  m_N-(n+1)\lambda,\,\hat h_x''\ge \bigl(a_Nn-1\bigr)\lambda\Bigr\}.
\end{equation}
Let~$A_n$ denote the event that the set in the union labelled by~$n$ is non empty and, given also~$\kappa>0$, abbreviate $B_n:=\{|\Gamma_N'(n\lambda)|>\texte^{\kappa n\lambda}\}$. By Proposition~\ref{prop-tight} there is~$\beta>0$ such that for all~$\kappa,\lambda$ large and all~$n\ge1$,  
\begin{equation}
P(B_n)\le\texte^{-\beta\kappa n\lambda}.
\end{equation}
Using the independence of~$h''$ of~$h'$ and the fact that $\Var(\hat h_x'') < C\textsl{}$, once $a_Nn>1$ a union bound and a standard Gaussian estimate show
\begin{equation}
P(A_n\setminus B_{n+1})\le \texte^{\kappa(n+1)\lambda}\texte^{-(a_Nn-1)^2\lambda^2/(2C)}.
\end{equation}
Putting these bounds together we get
\begin{equation}
P\Bigl(\,\Gamma_N(\lambda)\setminus\Gamma_N'(2\lambda)\ne\emptyset\Bigr)\le\sum_{n\ge2}\Bigl(\texte^{-\beta\kappa(n+1)\lambda}+ \texte^{\kappa(n+1)\lambda}\texte^{-(a_Nn-1)^2\lambda^2/(2C)}\Bigr).
\end{equation}
But~$a_N\to1$ as~$N\to\infty$, and so this sum converges uniformly as $N\to\infty$, and vanishes as $
\lambda\to\infty$, thus proving also the second part of the claim.
\end{proofsect}

Obviously, the extreme local maxima of~$h$ and~$h'$ will coincide only if the field~$\hat h''$ does not vary much in the neighborhood of these extreme points. This naturally leads us to study the oscillation of~$\hat h''$, where the oscillation of function~$g$ on a finite set~$\Lambda$ is defined by
\begin{equation}
\text{osc}_\Lambda\, g:=\max_{z\in\Lambda}g(z)-\min_{z\in\Lambda}g(z).
\end{equation}
The next claim to prove is thus:

\begin{lemma}
\label{lemma-T4}
For any $\lambda>0$, any~$\delta>0$ and any $r\ge1$,
\begin{equation}
\label{E:4.44}
\limsup_{N\to\infty}\,P\Bigl(\,\max_{x\in\Gamma'_N(\lambda)}\text{\rm osc}_{\Lambda_{2r}(x)}\hat h''>\delta\Bigr)=0.
\end{equation}
\end{lemma}

\begin{proofsect}{Proof}
Let us denote $F(x):=\{\max_{z\in\Lambda_{2r}(x)}|\hat h_z''-\hat h_x''|>\delta/2\}$. The event in \eqref{E:4.44} is then contained in $\bigcup_{x\in\Gamma_N'(\lambda)}F(x)$.
Given $\lambda > 0$, $M > 0$, $\rho \in (0,1)$ and $N \geq 1$, define
\begin{equation}
\label{E:4.46}
	A'_{N,M} := \bigl\{|\Gamma'_N(\lambda)|< M\bigr\}
	\quad\text{and} \quad
	B'_{N,\rho} := \bigl\{\Gamma'_N(\lambda) \subseteq V_{N, \rho}\bigr\} \,.
\end{equation}
Lemma~\ref{lemma-extremes-in-bulk} and Theorem~\ref{thm-DZ1} then imply
\begin{equation}
\label{E:4.46t}
\lim_{M \to \infty}\, \limsup_{N \to \infty}\, P(A^{\prime\,\cc}_{N, M}) = 0 
\quad \text{and} \quad
\lim_{\rho \downarrow 0}\, \limsup_{N \to \infty}\, P(B^{\prime\,\cc}_{N,\rho}) = 0 \,.
\end{equation}
On the other hand, by conditioning on $h'$ and using Chebyshev's inequality and a union bound, we get
\begin{equation}
P \Bigl(\,\,\bigcup_{x\in\Gamma_N'(\lambda)}F(x)\, \Big| \, A'_{N, M} \cap B'_{N, \rho} \Bigr) \,\leq\,
	4M |\Lambda_{2r}| \delta^{-2} \sup_{\substack{x,y \in V_{N, \rho}\\|x-y| \le 2r}}
			E\bigl(|\hat{h}''_x-\hat{h}''_y|^2\bigr).
\end{equation}
Lemma~\ref{lemma-var-bound} bounds the supremum by $C(\rho,r)/\log N$ and so the conditional probability tends to $0$ as $N \to \infty$ for any $\rho>0$ and any $M>0$. In light of \eqref{E:4.46t}, the claim follows.
\end{proofsect}

Consider now the a.s.\ well-defined mappings:
\begin{equation}
	\Pi(x):=\argmax_{\Lambda_{2r}(x)}h 
	\quad\text{and}\quad 
	\Pi'(x):=\argmax_{\Lambda_{2r}(x)}h' \,.
\end{equation}
Our next claim deals with the closeness of~$\Pi'(\Theta_{N,r})$ to~$\Theta_{N,r}'$, and $\Pi(\Theta_{N,r}')$ to $\Theta_{N,r}$, provided these are restricted to proper level sets:

\begin{lemma}
\label{lemma-T5}
The following holds with probability tending to one in the limits $N\to\infty$, $\delta\downarrow0$, $r\to\infty$ and~$\lambda\to\infty$ (in this given order):
\begin{equation}
\label{E:4.52}
	x\in\Theta_{N,r} \cap \Gamma_N(\lambda) \cap \Gamma_N'(\lambda)\quad\Rightarrow\quad\left\{
	\begin{aligned}&\Pi'(x)\in \Theta'_{N,r}, \quad
	\bigl|\Pi'(x)-x|\le \frac r2,\\
	&\text{and}\quad 0\le h_x-h_{\Pi'(x)}\le\delta,
	\end{aligned} \right.
\end{equation}
and
\begin{equation}
\label{E:4.53}
	x\in\Theta'_{N,r} \cap \Gamma_N(\lambda) \cap \Gamma_N'(\lambda)
	\quad\Rightarrow\quad\left\{
	\begin{aligned}&\Pi(x)\in\Theta_{N,r},\quad
	\bigl|\Pi(x)-x|\le \frac r2,\\
	&\text{and}\quad0\le h_{\Pi(x)} - h_x \le\delta.
	\end{aligned} \right.
\end{equation}
\end{lemma}

\begin{proofsect}{Proof}
For $N \geq 1$, $\lambda > 0$ and $r > 0$ abbreviate
\begin{equation}
\label{E:4.51}
	A_{N,\lambda,r}:=\Bigl\{\forall u,v\in \Gamma_N(\lambda)\colon |u-v|\le\tfrac r2 \ \text{\small\ or\ \  }|u-v|>{2r}\Bigr\} 
\end{equation}
and let $A'_{N,\lambda, r}$ be defined similarly with $\Gamma_N$ replaced by $\Gamma_N'$. By Theorem~\ref{thm-DZ2} we know that the probability of both of these events tends to one in the above limit and so we may henceforth assume that~$A_{N,\lambda,r}\cap A_{N,\lambda,r}'$ occurs.

Let~$x\in\Theta_{N,r}\cap\Gamma_N(\lambda)\cap\Gamma_N'(\lambda)$. We know that $h_x,h_x'\ge m_N-\lambda$ and so, in light of occurrence of $A_{N,r}\cap A_{N,r}'$, we have $h_u,h_u'<m_N-\lambda$ whenever $\ffrac r2<|u-x|\le 2r$. Since $\Pi'(x)$ is the point where~$h'$ is maximized on~$\Lambda_{2r}(x)$, we have $h'_{\Pi'(x)}\ge m_N-\lambda$ and so $\Pi'(x)\in\Lambda_{r/2}(x)$. It follows that $\Pi'(x)\in\Theta_{N,r}'$ as claimed. To get also the inequality between $h_x$ and~$h_{\Pi'(x)}$ we note
\begin{equation}
0\le h_x-h_{\Pi'(x)}=\hat h'_x-\hat h'_{\Pi'(x)}+\hat h''_x-\hat h''_{\Pi'(x)}\le\text{osc}_{\Lambda_{2r}}\hat h'',
\end{equation}
where for the second inequality we observed that $\hat h'_x-\hat h'_{\Pi'(x)}\le0$ and then applied the definition of oscillation. But Lemma~\ref{lemma-T4} guarantees that, with probability tending to one, $\text{osc}_{\Lambda_{2r}}\hat h''<\delta$ once~$N$ is large. This shows~\eqref{E:4.52}; the arguments for~\eqref{E:4.53} are completely analogous.\end{proofsect}

\begin{proofsect}{Proof of Proposition~\ref{lemma-T3}}
Fix $\epsilon > 0$ and let~$f$ be as given. Thanks to~$f$ having compact support, there is~$\lambda_0>0$ such that $f$ vanishes outside of $[0,1]^2 \times [-\lambda_0/2, \infty)$. Moreover, by uniform continuity of~$f$, for each~$M>0$ and~$r\ge1$ there is~$\delta>0$ such that, for all~$N$ large,
\begin{equation}
\label{E:4.53ww}
|x-x'|\le\frac r2\,\ \text{\ \small and\ }\ \, |h-h'|\le\delta\quad\Rightarrow\quad \Bigl|f\bigl(\tfrac xN,h\bigr)-f\bigl(\tfrac {x'}N,h'\bigr)\Bigr|<\epsilon/M.
\end{equation}
Let~$A_{N,M,r,\lambda,\delta}$ stand for the intersection of the event in \eqref{E:4.51} with the events $\{\Gamma_{N}(\lambda/2)\subseteq\Gamma_N'(\lambda)\}$ and $\{|\Gamma_N'(\lambda)|\le M\}$. We claim that once~$\lambda\ge\lambda_0$, on $A_{N,M,r,\lambda,\delta}$,
\begin{equation}
\label{E:4.56}
\begin{aligned}
\sum_{x\in\Theta_{N,r}} f\bigl(\tfrac xN,h_x-m_N\bigr) & = 
\sum_{x\in\Theta_{N,r}\cap\Gamma_N(\lambda)\cap\Gamma_N'(\lambda)}f\bigl(\tfrac xN,h_x-m_N\bigr) 
\\*[2mm]
& \le\frac\epsilon M\bigl|\Gamma_N'(\lambda)\bigr|+\sum_{x\in\Theta_{N,r}\cap\Gamma_N(\lambda)\cap\Gamma_N'(\lambda)}f\Bigl(\,\tfrac{\Pi'(x)}N,\,h_{\Pi'(x)}-m_N\Bigr)
\\*[2mm]
&\le\epsilon+\sum_{x\in\Theta'_{N,r}}f\bigl(\tfrac xN,h_x-m_N\bigr).
\end{aligned}
\end{equation}
Indeed, using $f\ge0$, the first inequality comes from inserting (freely, thanks to~$\lambda\ge\lambda_0$) the restriction $x\in\Gamma_N(\lambda/2)$ into the sum and then applying that~$\Gamma_N(\lambda/2)\subseteq\Gamma_N(\lambda)\cap\Gamma_N'(\lambda)$ on $A_{N,M,r,\lambda}$. In the second line we replace the arguments of~$f$ with the help of \eqref{E:4.53ww} and the estimates in \eqref{E:4.52}. Invoking the assumed bound $|\Gamma_N'(\lambda)|\le M$, the error is thus at most~$\epsilon$; the last inequality then follows from the first conclusion on the right of \eqref{E:4.52} and positivity of~$f$ again.

Defining, similarly, $A_{N,M,r,\lambda,\delta}'$ to be the intersection of the event in \eqref{E:4.51} with $\{\Gamma_{N}'(\lambda/2)\subseteq\Gamma_N(\lambda)\}$ and $\{|\Gamma_N(\lambda)|\le M\}$, we infer
\begin{equation}
\label{E:4.57}
\sum_{x\in\Theta'_{N,r}} f\bigl(\tfrac xN,h_x-m_N\bigr) 
\le\epsilon+\sum_{x\in\Theta_{N,r}}f\bigl(\tfrac xN,h_x-m_N\bigr),\qquad\text{on }A_{N,M,r,\lambda,\delta}'.
\end{equation}
Hence, for~$\lambda\ge\lambda_0$ the event in the claim is disjoint from~$A_{N,M,r,\lambda,\delta}\cap A_{N,M,r,\lambda,\delta}'$ and so to get the proof finished it suffices to have
\begin{equation}
	\lim_{\lambda \to \infty}\, \limsup_{r \to \infty}\, \limsup_{M \to \infty}\, \limsup_{\delta\downarrow0}\,\limsup_{N \to \infty}\,
		P\bigl( A_{N,M,r,\lambda,\delta}^\cc ) = 0,
\end{equation}
and similarly for the complement of $A_{N,M,r,\lambda,\delta}'$.
This follows by combining Lemma~\ref{lemma-T5} with Lemma~\ref{lemma-T2} and Theorem~\ref{thm-DZ1}.
\end{proofsect}

Before we move on to the proof of the distributional invariance, we note a simple fact concerning approximation of nearly-independent multivariate normals by independent ones:

\begin{lemma}
\label{lemma-cov}
Fix~$\sigma^2>0$ and let $Y_1,\dots,Y_n$ be i.i.d.\ $\NN(0,\sigma^2)$. For each continuous~$f\colon\R\to\R$ with compact support and each~$\epsilon>0$ there is~$\delta>0$ such that if $X=(X_1,\dots,X_n)$ are multivariate normal with $E(X_i)=0$ and 
\begin{equation}
\label{E:cov-ass}
\max_{1\le i,j\le n}\bigr|\text{\rm Cov}(X_i,X_j)-\sigma^2\delta_{ij}\bigl|<\delta
\end{equation}
then
\begin{equation}
\label{E:4.58ww}
\biggl|\log\frac{E\exp\{\sum_{i=1}^n f(X_i)\}}{\prod_{i=1}^n E\,\texte^{f(Y_i)}}\biggr|<\epsilon.
\end{equation}
\end{lemma}

\begin{proofsect}{Proof}
Writing $C$ for the covariance of~$(X_1,\dots,X_n)$, the fact that~$C$ is positive definite implies that $(X_1,\dots,X_n)=A(Y_1,\dots,Y_n)$ for an $n\times n$ matrix~$A$ such that $A^{\text{T}}A=\sigma^{-2}C$. In particular, as soon as \eqref{E:cov-ass} holds, $A$ is close to the identity once~$\delta$ is small enough. Plugging this into the numerator of \eqref{E:4.58ww} and applying uniform continuity of~$f$ along with tightness, the claim follows.
\end{proofsect}

\begin{proofsect}{Proof of Theorem~\ref{thm-1}}
Let $f$ be as stated. By the tightness of the maximum, we may and will suppose that
\begin{equation}
\label{E:4.61ww}
\max_{x\in V_N}|h_x'|,\,\max_{x\in V_N}|h_x''|\le 3\sqrt g\log N.
\end{equation}
Under such conditions $|h_x'-m_N|\le\lambda$ implies via Taylor expansion around $1$ of 
the square root in the definition of $\hat{h}'$ that
\begin{equation}
\label{E:4.62ww}
\biggl|\, h_x-\Bigl(h_x'-\frac\alpha2t+\sqrt{\frac t{g\log N}}\,h_x''\Bigr)\biggr|\le\frac{C}{\sqrt{\log N}}
\end{equation}
for some~$C=C(\lambda,t)>0$. This is similar to the calculation as in Section~\ref{sec2.1}. 

Consider the set $\Delta'_{N,\rho}(\lambda):=\{x\in V_{N,\rho}\colon |h_x'-m_N|\le\lambda\}$. Since~$f$ has compact support, the level sets of~$h$ and~$h'$ are related as in Lemma~\ref{lemma-T2}, and are confined into~$V_{N,\rho}$ by Lemma~\ref{lemma-extremes-in-bulk}. Using also the tightness of the maximum of $h'$, we thus have,
\begin{equation}
\lim_{\lambda\to\infty}\,\limsup_{r \to \infty}\,\limsup_{\rho\downarrow0}\,\limsup_{N\to\infty}\,P\biggl(\,\,
\sum_{x\in\Theta_{N,r}'\setminus\Delta'_{N,\rho}(\lambda)}\!\!f\bigl(\tfrac xN,h_x-m_N\bigr)>\epsilon\biggr)=0.
\end{equation}
This permits us to focus only on the sum of $f(\ffrac xN,h_x-m_N)$ over~$x\in\Theta_{N,r}'\cap \Delta'_{N,\rho}$.
As~$f$ is uniformly continuous, we now replace~$h_x$ by $h_x'-\frac\alpha2t+\hat h_x''$ and control the difference of the sums arising from \eqref{E:4.62ww} by the term $\delta_N|\{x\in V_N\colon f(\tfrac xN,h_x-m_N)>0\}|$, where~$\delta_N\to0$ as $N\to\infty$. As~$f$ has compact support, Theorem~\ref{thm-DZ1} ensures that the probability that this error exceeds a given~$\epsilon>0$ tends to zero as~$N\to\infty$. Using this in combination with Proposition~\ref{lemma-T3}, we conclude that
\begin{equation}
\label{E:4.40}
\lim_{\lambda\to\infty}\,\limsup_{r \to \infty}\,\limsup_{\rho\downarrow0}\,\limsup_{N\to\infty}\,P\biggl(\,\Bigl|\langle{\eta}_{N,r},f\rangle-\!\!\!\!
\sum_{x\in\Theta_{N,r}'\cap\Delta'_{N,\rho}(\lambda)}\!\!f\bigl(\tfrac xN,h_x'-m_N-\tfrac\alpha2t+\hat h_x''\bigr)\Bigr|>\epsilon\biggr)=0
\end{equation}
holds for any $\epsilon > 0$.

Our aim is to compute the limit of $E\texte^{-\langle{\eta}_{N,r},f\rangle}$. For this we replace $\langle{\eta}_{N,r},f\rangle$ by the sum in \eqref{E:4.40} and take the conditional expectation given~$h'$. For~$h'$ such that the event $A'_{N, \lambda, r}$ occurs, we have for all $u,v \in \Theta_{N,r}'\cap\Delta_{N,\rho}'(\lambda)$ with $u \neq v$,
\begin{equation}
	\text{Var}(\hat{h}_u'') = t + o(1) 
	\quad \text{and} \quad
	\text{Cov}(\hat h_u'',\,\hat h_v'')=o(1),\qquad N\to\infty.
\end{equation}
Therefore, by Lemma~\ref{lemma-cov}, on $A'_{N,\lambda, r} \cap \{|\Gamma'(\lambda)| < M\}$ we have
\begin{multline}
\qquad
E \biggl(\exp \Bigl\{-\sum_{x\in\Theta_{N,r}'\cap\Delta'_{N,\rho}(\lambda)} f\bigl(\tfrac xN,h_x'-m_N-\tfrac\alpha2t+\hat h_x''\bigr)\Bigr\} \,\bigg|\,h' \biggr) \\
\quad = \texte^{o(1)}\widetilde E\biggl( \exp \Bigl\{ -\sum_{x\in\Theta_{N,r}'\cap\Delta'_{N,\rho}(\lambda)}
 	f\bigl(\tfrac xN,h_x' - m_N + W^{(x)}_t -\tfrac{\alpha}{2}t \bigr)\Bigr\}\biggr)  
  \\
  = \exp \Bigl\{o(1)\,-\!\sum_{x\in\Theta_{N,r}'\cap\Delta'_{N,\rho}(\lambda)}
  		f_t\bigl(\tfrac xN, h_x'-m_N\bigr) \Bigr\},
\qquad N\to\infty,
\qquad
\end{multline}
where the expectation $\widetilde E$ is with respect to a measure under which $\{W_t^{(x)} \colon   
x\in\Theta_{N,r}'\cap\Delta'_{N,\rho}(\lambda)\}$ are i.i.d.\ Gaussians with mean $0$ and variance $t$, the function
$f_t$ is as defined in~\eqref{E:3.1} and the $o(1)$ term is random but tending to zero uniformly in~$h'$ under consideration.

Taking expectation also with respect to~$h'$ and noting that the function under expectation is at most one, Theorem~\ref{thm-DZ1}, Theorem~\ref{thm-DZ2} and the Bounded Convergence Theorem thus ensure
\begin{equation}
\lim_{\lambda\to\infty}\,\limsup_{r \to \infty}\,\limsup_{\rho\downarrow0}\,\limsup_{N\to\infty}\,\biggl|\,E\texte^{-\langle\eta_{N,r},f\rangle}-E\exp\Bigl\{-\!\!\!\!\!\sum_{x\in\Theta_{N,r}'\cap\Delta_{N,\rho}'(\lambda)}\!\!f_t\bigl(\tfrac xN,h_x'-m_N\bigr)\Bigr\}\biggr|=0 \,.
\end{equation}
We would like to drop the restriction to~$\Delta_{N,\rho}'(\lambda)$ and interpret the $N\to\infty$ limit (along the desired subsequence) of the second term using the limit process~$\eta$. Unfortunately, we cannot just roll the argument backward for~$f_t$ in place of~$f$ because~$f_t$ no longer has compact support. Notwithstanding, from the fact that~$f$ does, we get
\begin{equation}
f_t(x, h) \le C\texte^{-c h^2}.
\end{equation}
Therefore, on the event that 
$\{|\Gamma_N'(\theta)|<C\texte^{C\theta} \colon  \forall \theta>\lambda\} \, \cap \, \{h'_x \leq m_N + \lambda \colon \forall x \in V_N\}$, whose probability tends to $1$ as $N\to\infty$ and $\lambda \to \infty$ thanks to Proposition~\ref{prop-tight} and the tightness of the centered maximum, we have
\begin{equation}
\label{E:4.69ww}
\biggl|\,\sum_{x\in\Theta_{N,r}'\setminus\Delta'_{N,\rho}(\lambda)}f_t\bigl(\tfrac xN,h_x'-m_N\bigr)\biggr|
\le C\,\sum_{n\ge1}\texte^{C(\lambda+n)-c(\lambda+n-1)^2}.
\end{equation}
This is summable in $n$ and tends to $0$ as $\lambda \to \infty$. Using also Lemma~\ref{lemma-T1} we arrive at
\begin{equation}
\label{E:4.70}
\limsup_{N\to\infty}\,\bigl|\,
	E\texte^{-\langle\eta_{N,r_N}, f\rangle}- 
	E\texte^{-\langle\eta_{N,r_N}, f_t\rangle} \bigr|=0 \,.
\end{equation}
The bound \eqref{E:4.69ww} also implies that $f_t$ can be approximated by a continuous function $f_{t,\epsilon}$ with compact support so that $|\langle\eta_{N,r_N},f_t-f_{t,\epsilon}\rangle|\le\epsilon$.
Therefore, taking a limit of both terms in \eqref{E:4.70} along the subsequence for which $\eta_{N,r_N}$ converges to $\eta$, we recover~\eqref{E:3.2} as desired.
\end{proofsect}

\section{Multiple maxima and derivative martingale}
\label{sec5}
\noindent
At this stage of the proofs we know that all subsequential limits of the process $\eta_{N,r_N}$ take the form \eqref{E:1.9}; the last remaining issue is to establish uniqueness, i.e., that all such limits are actually the same. In Section~\ref{sec3.3} we have already reduced this to Theorem~\ref{thm-3} which deals with the joint law of maxima in several subsets of~$V_N$.
Here we prove  Theorem~\ref{thm-3} by a modification of the proof of Theorem 1.1 in Bramson, Ding and Zeitouni~\cite{BDingZ}.  As a by product, we will also prove Theorem~\ref{thm1.7}  giving~$Z(\textd x)$ the  interpretation of a derivative martingale. 

\subsection{Coupling to independent variables}
\label{sec5.1}\noindent
The argument underlying the proof of Theorem 1.1 in~\cite{BDingZ} relies on a coupling between the extreme process and a collection of i.i.d.\ random variables. This is based on the Gibbs-Markov property of the DGFF. We begin by a review of the relevant objects. 

Fix $K\in\N$ and recall our notation $w_i^K$, $i=1,\dots,K^2$, from Section~\ref{sec1.2} for points enumerating $\{x/K\colon x\in\{0,\dots,K-1\}^2\}$. If~$N$ is an natural number that is divisible by~$K$ and~$\delta>0$ is such that $\delta N/K\in\N$ we induce the sets
\begin{equation}
V_N^{K,i}:=N w_i^K+V_{N/K}\quad\text{and}\quad V_N^{K,\delta,i}:=N w_i^K+\bigl(N\delta/K,N(1-\delta)/K\bigr)^2\cap\Z^2
\end{equation}
Note that the boxes $V_N^{K,i}$ are disjoint and they tile~$V_N$ in the sense that  $V_N^{K,i}\cup\partial V_N^{K,i}$ cover $V_N\cup\partial V_N$. We also define $V_N^{K,\delta}:=\bigcup_{i=1}^{K^2}V_N^{K,\delta,i}$.

Passing to the continuum limit requires introducing also the continuum version of the above sets. In accord with~\cite{BDingZ} we recall the notation $B^{K,i} := w_i^K+(0,\ffrac1K)^2$ and set
\begin{equation}
	B^{K,\delta,i} := w_i^K+\bigl(\delta/K,(1-\delta)/K\bigr)^2\quad\text{and}\quad
	B^{K,\delta} := \bigcup_{i=1}^{K^2}B^{K,\delta,i}.
\end{equation}
These are the continuous analogues of the sets $V_N^{K,i}$, $V_N^{K,\delta,i}$ and $V_N^{K,\delta}$, respectively.
The reader is referred to Figure~1 in~\cite{BDingZ} for an illustration of a very similar geometric setup.

A key tool of \cite{BDingZ} is the representation of the DGFF as the sum $h = \hf + \hc$ of the ``fine field''~$\hf$ and the ``coarse field''~$\hc$ on $V_N$, which are independent of each other. Explicitly, set $\FF_{N,K} := \sigma \bigl(h_x \colon  x \in \bigcup_{i=1}^{K^2} \partial V_{N}^{K,i} \bigr)$ and define
\begin{equation}
	\hc := E \bigr(h \,\big| \FF_{N,K} \bigr) \quad\text{and}\quad
	\hf := h - \hc \,. 
\end{equation}
While the restrictions of $\hf$ to $V_N^{K,i}$, $i=1,\dots, K^2$ are independent copies of the  DGFF in $V_{N/K}$ with zero boundary condition on the boundary, $\hc$ admits a scaling limit
\begin{equation}
	\hc_{\lceil N \cdot \rceil}\, \underset{N\to\infty}\Lawlongarrow\, \Phi_K(\cdot),
\end{equation}
on $B^{K, \delta}$ for any $\delta > 0$, where $\bigl(\Phi_K(x) \colon  x \in B^{K,\delta} \bigr)$ is the Gaussian field defined in Section~\ref{sec1.2} and the convergence is with respect to the $L^\infty$-norm on bounded functions on $B^{K, \delta}$. (We will nonetheless use only the convergence of finite dimensional distributions.)

Abbreviate $h^\star:=\max_{x\in V_N}h_x$. The proof of the existence of the weak limit of the centered maximum $h^\star - m_N$ in \cite{BDingZ} is based on a comparison with an auxiliary process, defined in terms of~$\Phi_K$ and a collection of i.i.d.\ triplets $(\wp_i^K, Y_i^K, z_i^K)$ for $i=1, \dots, K^2$ as follows: Recall the probability density~$\psi$ on $[0,1]^2$ and constant~$C_\star$ from Theorem~\ref{thm-asymp}. For a sequence of numbers $b_K$ --- called $g(K)$ in \cite{BDingZ} --- to be determined later define $\wp_i^K \in \{0,1\}$, $Y_i^K \in (0,\infty)$ and $z_i^K \in B^{K,i}$ to be independent with marginal laws given, for $x \geq 0$ and Borel $A \subseteq [0,1]^2$, by
\begin{equation}
\label{E:5.4t}
\begin{aligned}
P(\wp_i^K  = 1) &:= C_\star \,b_K \,\texte^{-{\sqrt{2\pi}\, b_K}},\\
	P(Y_i^K \geq x) &:= \tfrac{b_K+x}{b_K}\texte^{-\sqrt{2\pi}\, x}
\end{aligned}
\end{equation}
and
\begin{equation}
\label{E:5.5t}
P\bigl(K(z_i^K -w^{K,i})\in A \bigr) := \int_A \psi(x) \,\textd x.
\end{equation}
Using these, we set 
\begin{equation}
G_i^{K,\delta} :=\begin{cases}
Y_i^K +b_K - 2\sqrt{g} \log K + \Phi_K(z_i^K),\qquad&\text{if }\wp_i^K = 1 \,\,\&\,\,z_i^K \in B^{K,\delta,i},
\\
-2\sqrt g\log K,\qquad&\text{otherwise}.
\end{cases}
\end{equation}
and denote $G_{K,\delta}^\star := \max_{i=1}^{K^2} G_i^{K,\delta}$. Let $\mu_N$ denote the law of $h^\star - m_{N}$ and $\mu_{K,\delta}$ stand for the law of~$G_{K,\delta}^\star$. Theorem~2.3 in \cite{BDingZ} then asserts that, for $b_K:=c\log\log K$ with some small~$c>0$, as $N\to\infty$, $K \to \infty$ and $\delta \downarrow 0$, the L\'evy distance between $\mu_N$ and $\mu_{K,\delta}$ tends to $0$. This yields weak convergence of $h^\star - m_N$ as the auxiliary process does not depend on $N$.
 

For our purposes, we need a generalization of Theorem~2.3 that addresses the joint law of maxima over several (different) subsets of~$[0,1]^2$. To this end, for any non-empty open~$A\subset[0,1]^2$ we first define
$h^\star_A$ to be the maximum of the GFF in $\{x \in V_N \colon  \tfrac{x}{N} \in A\}$. Then, given a collection $\underline{A} = (A_1, \dots, A_m)$ of  non-empty open subsets of $[0,1]^2$, we use $\mu_{N, \underline{A}}$ to denote the joint law of $\bigl(h^\star_{A_l} - m_N \bigr)_{l=1}^m$. For the above auxilary process we similarly define 
\begin{equation}
	G_{K, \delta, A}^\star := \max \bigl\{G_i^{K,\delta} \colon  B^{K,i} \subseteq A\bigr\}
\end{equation}
and write $\mu_{K,\delta, \underline{A}}$ for the joint law of $\bigl(G_{K, \delta, A_l}^\star\bigr)_{l=1}^m$. Our aim is then to prove:

\begin{theorem}
\label{our-thm2.3}
Set $b_K:=c\log\log K$ for some~$c>0$ sufficiently small.
For any $m\ge1$, let $d^{(m)}$ denote the L\'{e}vy metric on probability measures on $\R^m$.
Then for any $m\ge1$ and any collection $\underline{A} = (A_1, \dots, A_m)$ of non-empty open subsets of $[0,1]^2$,
\begin{equation}
\label{E:5.8q}
	\lim_{\delta \downarrow 0}\, \limsup_{K \to \infty}\, \limsup_{N \to \infty} \ 
		d^{(m)} \bigl(\mu_{N, \,\underline{A}}, \, \mu_{K, \delta,\,\underline{A}} \bigr) = 0.
\end{equation} 
\end{theorem}

Before we set out to prove this, let us show how this gets us the convergence of the joint law of maxima in a collection of open sets:

\begin{proofsect}{Proof of Theorem~\ref{thm-3}}
From \eqref{E:5.8q} and the fact that $\mu_{K, \delta,\,\underline{A}}$ does not depend on~$N$ we infer that $\mu_{N, \underline{A}}$ is Cauchy in the L\'evy metric $d^{(m)}$. In particular, the joint law of $\bigl(h^\star_{A_l} - m_N \bigr)_{l=1}^m$ on~$\R^m$ converges weakly as $N \to \infty$ for any non-empty open  $A_1,\dots,A_m\subset\R$.
\end{proofsect}

Similar arguments yield also the characterization of~$Z(\textd x)$ as a derivative martingale:

\begin{proofsect}{Proof of Theorem~\ref{thm1.7}}
Let~$Z(\textd x)$ be the random measure characterizing the limiting point process and let $Z_{K,\delta}$ denote the measure from \eqref{E:1.12t} with $\epsilon_K$ replaced by~$\delta$. 

By standard approximation arguments, in order to prove $Z_{K,\delta}\Lawarrow Z$, it suffices to show that the joint law of $(Z_{K,\delta}(A_1),\dots,Z_{K,\delta}K(A_m))$ on~$\R^m$ converges weakly to that of $(Z(A_1),\dots,Z(A_m))$, for any collection of disjoint non-empty open sets $A_1,\dots,A_m$ with $Z(\partial A_i)=0$ a.s.\ for all $i=1,\dots,m$ and such that $\dist(A_i,A_j)>0$ for any $i\ne j$. Under such conditions, \eqref{E:3.27} and Theorem~\ref{our-thm2.3} yield
\begin{equation}
\begin{aligned}
\label{E:5.9r}
E \biggl(\exp \Bigl\{ -{\alpha}^{-1} \sum_{l=1}^m \texte^{-\alpha t_l} Z(A_l) \Bigr\}\biggr) 
&=\lim_{\delta\downarrow0}\lim_{K\to\infty} E\Bigl(\,\,\prod_{l=1}^m\,\prod_{i\colon B^{K,i}\subset A_l}\1_{\{G_i^{K,\delta}\le t_l\}}\Bigr)
\\
&=\lim_{\delta\downarrow0}\lim_{K\to\infty} E\Bigl(\,\,\prod_{l=1}^m\,\prod_{i\colon B^{K,i}\cap A_l\ne\emptyset}\1_{\{G_i^{K,\delta}\le t_l\}}\Bigr),
\end{aligned}
\end{equation}
where the second line follows by taking the first line for~$A_i$ enlarged by some~$\epsilon>0$ such that $\epsilon<\min_{i\ne j}\dist(A_i,A_j)$ and then taking~$\epsilon\downarrow0$ and applying that $Z(\partial A_i)=0$ a.s.\ for all~$i$. Our strategy is to represent the right-hand side using the measures~$Z_{K,\delta}$.

We first condition on~$\Phi_K$ and examine the conditional distribution function of~$G_i^{K,\delta}$. Here we note that if $K$ is so large that $t>-2\sqrt g\log K$
then $G_i^{K,\delta}>t$ is equivalent to $\wp_i^K=1$, $z_i^K\in B^{K,\delta,i}$ and $Y_i^K +b_K - 2\sqrt{g} \log K + \Phi_K(z_i^K)> t$. Moreover, since $Y_i^K>0$ a.s., the last condition is true trivially when $b_K - 2\sqrt{g} \log K + \Phi_K(z_i^K)\ge t$. Invoking \twoeqref{E:5.4t}{E:5.5t}, we get
\begin{multline}
\label{E:5.10u}
\qquad
P\bigl(G_i^{K,\delta}> t\big|\Phi_K)
\\=C_\star\int_{(\delta,1-\delta)^2}\!\!\textd z\,\psi(z)\,\widetilde F_K\Bigl(\bigl(t+2\sqrt{g} \log K - \Phi_K(w^{K,i}+z/K)\bigr)\vee b_K\Bigr),
\qquad
\end{multline}
where $\widetilde F_K(s):=s\texte^{-\alpha s}$ (with $\alpha:=2/\sqrt g =\sqrt{2\pi}$ as usual).

Since the argument of $\widetilde F_K$ is at least~$b_K$, the integral \eqref{E:5.10u} is at most $\beta_K:=C_\star b_K\texte^{-\alpha b_K}$. Note that $\beta_K\to0$ once $b_K\to\infty$. Using that $\{G_i^{K,\delta}\colon i=1,\dots,K^2\}$ are independent conditional on~$\Phi_K$, standard comparisons of $1-x$ with $\texte^{-x}$ for small~$x$ and a monotonicity argument yield
\begin{equation}
\label{E:5.12r}
\lim_{K\to\infty}E\Bigl(\,\,\prod_{l=1}^m\,\prod_{i\colon B^{K,i}\subset A_l}\1_{\{G_i^{K,\delta}\le t_l\}}\Bigr)
=
\lim_{K\to\infty}
E\biggl(\exp\Bigl\{-\sum_{l=1}^m\sum_{i\colon B^{K,i}\subset A_l}P\bigl(G_i^{K,\delta}> t_l\big|\Phi_K\bigr)\Bigr\}\biggr),
\end{equation}
and similarly for $B^{K,i}\subset A$ replaced by $B^{K,i}\cap A\ne\emptyset$. In light of \eqref{E:5.9r}, to prove the claim it thus suffices to show that, for any non-empty open set $A\subset[0,1]^2$ with $Z(\partial A)=0$ a.s.\ and any $t\in\R$,
\begin{equation}
\label{E:5.13w}
\alpha^{-1}\texte^{-\alpha t}\sum_{i\colon B^{K,i}\subset A} Z_{K,\delta}\bigl(B^{K,i}\bigr)-\sum_{i\colon B^{K,i}\subset A}P\bigl(G_i^{K,\delta}> t\big|\Phi_K)
\,\overset{P}{\underset{\begin{subarray}{c}
K\to\infty\\\delta\downarrow0
\end{subarray}}
\longrightarrow}\,0
\end{equation}
and similarly for $B^{K,i}\subset A$ replaced by $B^{K,i}\cap A\ne\emptyset$. (Note that \eqref{E:5.9r} and \eqref{E:5.12r} then ensure that the difference between the two kinds of sums tends to zero in probability.)

For $F(s)$ 
as in the statement of the theorem, once $\Phi_K(w^{K,i}+z/K)\le2\sqrt g\log K-b_K+t$ holds, we get
\begin{multline}
\qquad
\widetilde F_K\Bigl(\bigl(t+2\sqrt{g} \log K - \Phi_K(w^{K,i}+z/K)\bigr)\vee b_K\Bigr)
\\=\texte^{-\alpha t} \bigl(1+O(t/b_K)\bigr) F\bigl(2\sqrt{g} \log K - \Phi_K(w^{K,i}+z/K)\bigr).
\qquad
\end{multline}
Setting $c_\star:=\alpha C_\star$, comparing \eqref{E:5.10u} with the sum in \eqref{E:5.13w} (as well as their alternatives with $B^{K,i}\cap A\ne\emptyset$) and using that both $\widetilde F_K$ and $F_K$ are bounded, it is enough to show that
\begin{equation}
\label{E:5.14}
\sum_{i=1}^{K^2}\int_{(\delta,1-\delta)^2}\!\!\textd z\,\psi(z)\1_{\{\Phi_K(w^{K,i}+z/K)\ge2\sqrt g\log K-b_K+t\}}\overset{P}{\underset{K\to\infty}\longrightarrow}\,0
\end{equation}
for any $t\in\R$. Noting that $\text{Var}(\Phi_K(w^{K,i}+z/K))=g\log K+O(1)$ while $b_K^2=o(\log K)$, the standard Gaussian tail estimate yields
\begin{equation}
P\bigl(\Phi_K(w^{K,i}+z/K)\ge2\sqrt g\log K-b_K+t\bigr)\le \frac{c_1\texte^{\alpha b_K}}{K^2\sqrt{\log K}}
\end{equation}
for some $c_1=c_1(t)$. But $b_K$ was so far arbitrary subject to $b_K\to\infty$ and $b_K\le c\log\log K$ for~$c$ small enough, so we can make the choice such that the right-hand side is $o(K^2)$ as $K\to\infty$. Taking expectation with respect to~$\Phi_K$ we hereby get \eqref{E:5.14} and thus the resulting claim.
\end{proofsect}

\subsection{Proof of Theorem~\ref{our-thm2.3}}
The proof of Theorem 2.3 in~\cite{BDingZ} relies on further lemmas and propositions, so we first need to produce proper modifications for them.  The first of these is a tightness result for $h^\star_A - m_N$ in analogy to the tightness of $h^\star - m_N$.

\begin{lemma}
The sequence $h^\star_{A} - m_N$ is tight.
\end{lemma}

\begin{proofsect}{Proof}
Since $A$ is open and non-empty, there exists $K$ and $i$ such that $B^{K,i} \subseteq A$. Fixing $\delta \in (0,1)$ and letting $z_i \in V_N^{K,\delta, i}$ be such that $\hf_{z_i} = \max_{v \in V_N^{K,\delta, i}}\, \hf$, we have
\begin{equation}
	\hf_{z_i} + \hc_{z_i} \leq h^\star_{A} \leq h^\star \,.
\end{equation}
The tightness of $h^\star$ (proved in~\cite{BDeuG} and~\cite{BrZ}) implies that 
\begin{equation}
\lim_{t \to \infty} \,\limsup_{N \to \infty} \, P \bigl( h^\star_{A} - m_N > t \bigr) 
	\leq \lim_{t \to \infty}\, \limsup_{N \to \infty}\, P \bigl( h^\star - m_N > t \bigr) = 0
\end{equation}
giving the tightness for the upper tail of $h^\star_{A}$.

Concerning the lower tail we note that, since $\hf$ restricted to $V_N^{K,i}$ is a GFF on a translate of $V_{N'/K}$ and due to Lemma~\ref{lemma-extremes-in-bulk} and Theorem~\ref{thm2.1}, we have
\begin{equation}
\label{E:1.3}
	\lim_{\delta\downarrow0} \,\limsup_{t \to \infty}\, \limsup_{N \to \infty} \ 
	P \bigl( \hf_{z_i} -m_{N/K} < -t \bigr) = 0.
\end{equation}
At the same time, for all $x \in V_N^{K,\delta, i}$ we have
\begin{equation}
\begin{split}
	\text{Var}(\hc_{x}) = & \text{Var}\bigl(E(h_x | \FF_{N,K}) \bigr)
	= \text{Var}\bigl(h_x\bigr) - E\bigl(\text{Var}(h_x | \FF_{N,K})\bigr)  \\
	& \leq g \log N - g \log (N'/K) + C \leq g \log K + C' \,,
\end{split}
\end{equation}
where $C'$ depends on $\delta$.  This follows from Lemma~\ref{lemma-var-bound} (of the present paper) since $x$ is away from the boundary of both $V_N$ and $ V_N^{K,\delta, i}$. Conditioning on $z_i$ and using the independence of $\hf$ and $\hc$ we have
\begin{equation}
\label{E:1.4}
	P \bigl (\hc_{z_i} < -t \bigr) \leq C \exp \Bigl(-\tfrac{t^2}{2(g \log K + C')} \Bigr) 
\end{equation}
uniformly in $N$. Writing
\begin{equation}
	P \bigl( h^\star_{A} - m_N < -t \bigr)
		\leq P \bigl( \hf_{z_i} - m_{N/K} < -t/2 + 2\sqrt{g} \log K \big) +
			 P \bigl( \hc_{z_i} < -t/2 \bigr)
\end{equation}
and using~\eqref{E:1.3} and~\eqref{E:1.4}, the result follows. 
\end{proofsect}
	 	
Next we need proper modifications of Propositions~5.1 and~5.2 and Lemma~6.4 in~\cite{BDingZ}. Given a non-empty open $A\subset [0,1]^2$, we define
$V_{N,A}^{K,\delta} := \bigcup \{V_N^{K,\delta, i} \colon  B^{K,i} \subseteq A \}$ and
$\Delta_{N,A}^{K,\delta} := \{x \in V_N \colon  x/N \in A \} \setminus V_{N,A}^{K,\delta}$. Notice that since $A$ is open, the Lebesgue measure of $A \setminus \bigcup\{B^{K,i} \colon  B^{K,i} \subseteq A\}$ tends to $0$ as $K \to \infty$. It follows that
\begin{equation}
	\lim_{\delta\downarrow0}\, \limsup_{K \to \infty}
	\limsup_{N \to \infty}
	\, \frac{|\Delta_{N,A}^{K,\delta}|}{|V_N|} = 0.
\end{equation}
Our substitute for Proposition~5.1 in~\cite{BDingZ}, which shows that the maximum of~$h$ is unlikely to fall inside the set $\Delta_N^{K,\delta}$, is then:

\begin{proposition}
\label{prop-5.1'}
Let $A$ be any open non-empty subset of $[0,1]^2$. Then
\begin{equation}
	\lim_{\delta\downarrow0} \,\limsup_{K \to \infty}\, \limsup_{N \to \infty} \ 
		P \Bigr(\max_{x \in V_{N,A}^{K,\delta}}h_x \neq h^\star_A \Bigr) = 0 \,.
\end{equation}
\end{proposition}

\begin{proofsect}{Proof}
The proof is the same, only that we use $\Delta_{N,A}^{K,\delta}$ instead of $\Delta_N$ and 
the tightness of $h^\star_{A}-m_N$ instead of that of $h^\star$.
\end{proofsect}

In place of Proposition 5.2 of~\cite{BDingZ}, which shows that the maximum of~$h$ occurs with high probability at the maximum of the fine field, we in turn use:

\begin{proposition}
\label{prop:5.2'}
Let $A \subset [0,1]^2$ be non-empty and open. Let $z_i$ be such that $\max_{x \in V_N^{K,\delta, i}} \hf_x = h_{z_i}$ and $\bar{z}$ be such that $\max \{ h_{z_i} \colon  B^{K,i} \subseteq A\} = h_{\bar{z}}$. Then for any $\epsilon, \delta > 0$ fixed
\begin{equation}
	\lim_{K \to \infty} \limsup_{N \to \infty} \ 
		P \Bigl( \max_{x \in V_{N,A}^{K,\delta}} h_x \geq h_{\bar{z}} + \epsilon \Bigr) = 0 \,.
\end{equation}
Furthermore, there exists a sequence $(b_K)_{K \geq 1}$ with $b_K \to \infty$ as $K \to \infty$ such that
\begin{equation}
	\lim_{K \to \infty} \limsup_{N \to \infty} \ 
		P \bigl( \hf_{\bar{z}} \leq m_{N/K} + b_K \bigr) = 0 \,.
\end{equation}
\end{proposition}

\begin{proofsect}{Proof}
The proof is the same, only that we use the tightness of $h^\star_{A}-m_N$ (and Proposition~\ref{prop-5.1'}) 
instead of that of $h^\star$.
\end{proofsect}

Finally, we need to adapt Lemma 6.4 of~\cite{BDingZ} which shows that small changes in the coordinate of the coarse field do not affect the resulting maximum. For our purposes, this becomes:

\begin{lemma}
\label{lemma:6.4'}
For $i=1, \dots, K^2$, let $z_i$ be as in Proposition~\ref{prop:5.2'} and $z'_i \in V_N^{K,\delta,i}$ be measurable with respect to $\{\hf_x \colon  x \in V_N^{K,i}\}$ such that $|z_i - z'_i| = o(N/K)$ as $N,K \to \infty$. Then for any $A \subseteq [0,1]^2$ open and non-empty
\begin{equation}
	\lim_{K \to \infty} \limsup_{N \to \infty} \ 
		d \Bigl( \max \bigr\{ \hf_{z_i} + \hc_{z_i} \colon  B^{K,i} \subseteq A \bigr\} ,\,
			     \max \bigr\{ \hf_{z_i} + \hc_{z'_i}\colon  B^{K,i} \subseteq A \bigr\} \Bigr) = 0 \,.		
\end{equation}
where $d$ is the L\'{e}vy metric for $\R$-valued random variables.
\end{lemma}

\begin{proofsect}{Proof}
The proof of Lemma~6.4 from \cite{BDingZ} can be used without significant changes.
\end{proofsect}

We are now ready to prove the main theorem of this section:

\begin{proofsect}{Proof of Theorem~\ref{our-thm2.3}}
We follow the original proof of Theorem 2.3. For each $A_1, \dots, A_m$, define
\begin{equation}
	\bar h_{A_l}^\star := \max \bigl\{ \hf_{z_i} + \hc_{z_i} \colon  
		B^{K,i} \subseteq A,\, \hf_{z_i} > m_{N/K} + b_K \bigr\} \,.		
\end{equation}
Then for all $l=1, \dots, m$, we have by Proposition~\ref{prop-5.1'} and Proposition~\ref{prop:5.2'},
\begin{equation}
	P \bigl( \bigl| \bar{h}^\star_{A_l} - h^\star_{A_l} \bigr| > \epsilon \bigr) < \epsilon \,.
\end{equation}
Denoting by $\nu_{N, \underline{A}}^{K,\delta}$ the law of $\bigl(\bar{h}_{A_l}^\star\bigr)_{l=1}^m$, this shows $d \bigl( \mu_{N,\underline{A}} \,,\,\, \nu_{N,\underline{A}}^{K,\delta} \bigr) < \epsilon$. 

Next we employ the same coupling of $h$ with $(\wp_i^K, Y_i^K, z_i^K)_{i=1}^{K^2}$, as in the original proof. Proposition 6.3 from~\cite{BDingZ} together with Lemma~\ref{lemma:6.4'} imply that the law of 
\begin{equation}
	\Bigl( \max \bigl\{ b_K + Y_i^K + \hc_{z_i^K} - (m_N - m_{N/K}) \colon   
		\wp_i^K = 1 ,\, z_i^K \in V_N^{K,\delta,i} ,\, B^{K,i}_N \subseteq A_l \bigr\} \Bigr)_{l=1}^m,
\end{equation}
to be denoted by $\bar{\nu}_{N,A}^{K,\delta}$, satisfies 
$d^{(m)} \bigl( \bar{\nu}_{N,\underline{A}}^{K,\delta} \,,\,\, \nu_{N,\underline{A}}^{K,\delta} \bigr) < \epsilon$, for all $K$,$N$ large enough. Using the convergence of $\hc$ to $\Phi_K$ as $N \to \infty$ and triangle inequality, we conclude that
\begin{equation}
	\lim_{\delta\downarrow0}\, \limsup_{K \to \infty}\, \limsup_{N \to \infty} \,  
	 d^{(m)} \bigl( \mu_{N, \underline{A}} \, , \,\, \mu_{K,\delta, \underline{A}} \bigr) = 0 \,,
\end{equation}
as desired. Finally, since $\mu_{K,\delta, \underline{A}}$ does not depend on~$N$, this implies that $\mu_{N, \underline{A}}$ is Cauchy and hence converges weakly to a limit. 
\end{proofsect}

\section*{Acknowledgments}
\noindent
This research has  been partially supported by NSF grant DMS-1106850, NSA grant H98230-11-1-0171 and GA\v CR projects P201/11/1558 and P201/12/2613. The authors wish to thank J.~Ding, T.~Liggett and O.~Zeitouni for helpful discussions.

\end{document}